\documentclass[11pt,dvips]{article}
\usepackage{latexsym}
\usepackage{amsmath,amsthm,amsfonts,amssymb}
\usepackage[T1]{fontenc}
\usepackage[latin1]{inputenc}
\usepackage{aeguill}
\usepackage{url}
\usepackage{enumerate}
\usepackage{mdwlist} 
\usepackage{rotating} 
\usepackage{graphicx}
\usepackage[a4paper,textwidth=15cm,textheight=25cm]{geometry}
\usepackage{export}

\newtheorem{thm}{Theorem}[section]

\newtheorem*{thm*}{Theorem}
\newtheorem{dfn}[thm]{Definition} 
\newtheorem*{dfn*}{Definition}

\newtheorem{cor}[thm]{Corollary}
\newtheorem*{cor*}{Corollary}

\newtheorem{prop}[thm]{Proposition} 
\newtheorem*{prop*}{Proposition} 
\newtheorem*{properties*}{Properties} 
 
\newtheorem{lem}[thm]{Lemma} 
\newtheorem*{lem*}{Lemma}

\newtheorem*{claim*}{Claim} 
 
\newtheorem*{fact*}{Fact} 
\newtheorem{fact}[thm]{Fact}

\theoremstyle{remark}
\newtheorem*{rem*}{Remark}
\newtheorem{rem}[thm]{Remark}
\newtheorem*{example*}{Example}

\newenvironment{SauveCompteurs}[1]{%
\newcommand{\monparametre}{#1}
\openexport{\monparametre_sauve}
  \Export{thm}\Export{section}\Export{subsection}\Export{subsubsection}
\closeexport}{}

\newenvironment{UtiliseCompteurs}[1]{%
\newcommand{\monparametre}{#1}
\openexport{\monparametre_aux}
  \Export{thm}\Export{section}\Export{subsection}\Export{subsubsection}
\closeexport
\PackageInfo{export}{\MessageBreak
Importations from \monparametre_sauve.xpt\MessageBreak}%
\InputIfFileExists{\monparametre_sauve.xpt}{\relax}{\relax}%
\renewcommand{\label}[1]{}
}{%
\PackageInfo{export}{\MessageBreak
Importations from \monparametre_aux.xpt\MessageBreak}%
\InputIfFileExists{\monparametre_aux.xpt}{\relax}{\relax}}

\newlength{\espaceavantspecialthm}
\newlength{\espaceapresspecialthm}
{\normalfont \vskip \espaceapresspecialthm}

\newenvironment{specialthm*}[1]{
\vskip\espaceavantspecialthm \noindent \textbf{#1} \itshape}%
{\normalfont \vskip \espaceapresspecialthm}

\newlength{\espaceavantenonce}
\newlength{\espaceapresenonce}
\newcommand{\fontetitreun}[1]{\textbf{#1}} 
\newcommand{\fontetitredeux}[1]{\textit{#1}} 

{\normalfont \vskip \espaceapresenonce}

\newenvironment{enonce1*}[1]{
\vskip\espaceavantenonce \noindent \fontetitreun{#1} \itshape}%
{\normalfont \vskip \espaceapresenonce}

{\vskip \espaceapresenonce}

\newenvironment{enonce2*}[1]{
\vskip\espaceavantenonce \noindent \fontetitredeux{#1} }%
{\vskip \espaceapresenonce}

\makeatletter
\edef\@tempa#1#2{\def#1{\mathaccent\string"\noexpand\accentclass@#2 }}
\@tempa\rond{017}
\makeatother

\newcommand{\es}{\emptyset}
\renewcommand{\phi}{\varphi} 
\newcommand{\m} {^{-1}} 
\newcommand{\eps} {\varepsilon}

\newcommand {\ra} {\rightarrow}

\newcommand {\onto} {\twoheadrightarrow}

\newcommand{\actson}{\,\raisebox{1.8ex}[0pt][0pt]{\begin{turn}{-90}\ensuremath{\circlearrowright}\end{turn}}\,}

\newcommand{\ol}[1]{\overline{#1}}

\newcommand{\normal} {\vartriangleleft}
\renewcommand{\subsetneq}{\varsubsetneq}

\newcommand{\dunion}{\sqcup}

\newcommand{\ie} {i.~e.\ }

\newcommand {\cale} {{\mathcal {E}}}   
\newcommand {\calf} {{\mathcal {F}}}   
\newcommand {\calg} {{\mathcal {G}}}   
\newcommand {\calh} {{\mathcal {H}}}

\newcommand {\calm} {{\mathcal {M}}}

\newcommand {\cals} {{\mathcal {S}}}   
\newcommand {\calt} {{\mathcal {T}}}

\newcommand {\caly} {{\mathcal {Y}}}

\newcommand {\bbR} {{\mathbb {R}}}

\newcommand {\bbZ} {{\mathbb {Z}}}

\newcommand{\dom} {\mathop{\mathrm{dom}}}
\newcommand{\Isom} {\mathop{\mathrm{Isom}}}
\renewcommand{\Im} {\mathop{\mathrm{Im}}}

\newcommand{\id} {\mathrm{id}}

\usepackage[all]{xypic}

\newcommand{\Xt}{X^{(t)}}

\newcommand{\Tmin}{T_\mathrm{min}}

\begin{document}

\title{Actions of finitely generated groups on $\bbR$-trees.}
\author{Vincent Guirardel}

\maketitle

\begin{abstract}
We study actions of finitely generated groups on $\bbR$-trees under some stability hypotheses.
We prove that either the group splits over some controlled subgroup (fixing an arc in particular), 
or the action can be obtained by gluing together actions of simple types:
actions on simplicial trees, actions on lines, and actions coming from measured foliations on $2$-orbifolds.
This extends results by Sela and Rips-Sela. 
However, their results are misstated, and we give a counterexample to their statements.

The proof relies on an extended version of Scott's Lemma of independent interest. 
This statement claims that if a group $G$ is a direct limit of groups having suitably compatible splittings,
then $G$ splits.
\end{abstract}

Actions of groups on $\bbR$-trees are an important tool in geometric group theory.
For instance, actions on $\bbR$-trees are used to compactify sets of 
hyperbolic structures (\cite{MS_valuationsI,Pau_topologie}) or 
Culler-Vogtmann's Outer space (\cite{CuMo,BF_outer,CoLu_very}).
Actions on $\bbR$-trees are also a main ingredient in Sela' approach to 
acylindrical accessibility \cite{Sela_acylindrical} (see \cite{Delzant_accessibilite} or \cite{KaWe_acylindrical} for alternative approaches)
and in some studies of morphisms of a given group to a (relatively) hyperbolic group \cite{Sela_hopf,BeleSzcz_endomorphisms,DrutuSapir_tree-graded}.
Limit groups and limits of (relatively) hyperbolic groups \cite{Sela_diophantine1,Sela_diophantine7,Gui_limit,Alibegovic_MR,Groves_limit_hypII},
and Sela's approach to Tarski's problem (\cite{Sela_diophantine6}) are also studied using $\bbR$-trees as a basic tool. 

To use $\bbR$-trees as a tool, one needs to understand the structure of a group acting on an $\bbR$-tree.
The main breakthrough in this analysis is due to Rips.
He proved that any finitely generated group acting freely on an $\bbR$-tree
is a free product of surface groups and of free abelian groups (see \cite{GLP1,BF_stable}).

More general results involve \emph{stability} hypotheses. 
Roughly speaking, these hypotheses 
prohibit infinite sequences of nested arc stabilizers.
A simple version of these stability hypotheses is the \emph{ascending chain condition}:
for any decreasing sequence of arcs $I_1\supset I_2\supset ...$ whose lengths converge to $0$, 
the sequence of their pointwise arc stabilizers $G(I_1) \subset G(I_2)\subset\dots$
stabilizes. The ascending chain condition implies BF-stability used below.
See Definition \ref{dfn_stable} for other versions of stability and relations between them.

Another theorem by Rips claims that if $G$ is finitely presented and has a BF-stable action on an $\bbR$-tree $T$,
then either $T$ has an invariant line or $G$ splits over an
extension of a subgroup fixing an arc by a cyclic group (see \cite{BF_stable,Gui_approximation}).
Using methods in the spirit of \cite{GLP1}, the author proved under the same hypotheses
that one can \emph{approximate}  the action $G\actson T$
by a sequence of actions on simplicial trees $G\actson T_k$ converging to $T$ (\cite{Gui_approximation}).
Edge stabilizers of $T_k$
are extensions by  finitely generated free abelian groups of  subgroups fixing an arc in $T$.
The convergence is in the length functions topology.
\\


Next results  give a very precise structure of the action on the $\bbR$-tree.
This is the basis for Sela's shortening argument (see for instance \cite{RiSe_structure,Sela_acylindrical}).
These results say that the action on the $\bbR$-tree can be obtained by gluing together simple building blocks.
The building blocks are actions on simplicial trees, 
actions on a line, and actions coming from a measured foliation on a $2$-orbifold.
Two building blocks are glued together along one point, and globally, the 
combinatorics of the gluing is described by a simplicial tree (see section \ref{sec_goa} for more details).
When $G\actson T$ arises in this fashion, we say that $T$ splits as a \emph{graph of actions on $\bbR$-trees}.

In \cite{RiSe_structure} and \cite{Sela_acylindrical}, Rips-Sela and Sela give a structure theorem for actions 
of finitely presented groups and of finitely generated groups on $\bbR$-trees.
The proof uses the following super-stability hypothesis:
an arc $I\subset T$ is \emph{unstable} if there exists $J\subset I$ with $G(I)\subsetneq G(J)$;
an action $G\actson T$ is \emph{super-stable} if the stabilizer of any unstable arc is trivial.
In particular, super-stability implies that chains of arc stabilizers have length at most 2.
Sela's result is as follows: consider a minimal action $G\actson T$ of a finitely generated group on an $\bbR$-tree,
and assume that this action is super-stable, and has trivial tripod stabilizers.
Then either $G$ splits as a free product,
or $T$ can be obtained from simple building blocks as above (\cite[Theorem 3.1]{Sela_acylindrical}).

Actually, Sela's result is stated under the ascending chain condition instead of super-stability.
However, the proof really uses the stronger hypothesis of super-stability. 
In section \ref{sec_example}, we give a counterexample to the more general statement.
This does not affect the rest of Sela's papers since super-stability is satisfied in the cases considered.
Since the counterexample is an action of a finitely presented group, this also shows that
one should include super-stability in the hypotheses of Rips-Sela's statement \cite[Theorem 10.8]{RiSe_structure}
(see also \cite[Theorem 2.3]{Sela_acylindrical}).
\\

Our main result generalizes Sela's result in two ways.
First, we don't assume that tripod stabilizers are trivial.
Moreover, we replace super-stability by the ascending chain condition together with weaker assumptions on stabilizers of unstable arcs.


 

\begin{thm*}[Main Theorem]
Consider a non-trivial minimal action of a finitely generated group $G$ on an $\bbR$-tree $T$ by isometries.
Assume that
\begin{enumerate*}
\item $T$ satisfies the ascending chain condition;
\item for any unstable arc $J$,
\begin{enumerate*}
\item $G(J)$ is finitely generated;
\item $G(J)$ is not a proper subgroup of any conjugate of itself
\ie $\forall g\in G$, $G(J)^g\subset G(J)\Rightarrow G(J)^g= G(J)$.
\end{enumerate*}
\end{enumerate*}

Then either $G$ splits over the stabilizer of an unstable arc or over the stabilizer of an infinite tripod,
or $T$ has a decomposition into a graph of actions where each
vertex action is either 
\begin{enumerate*}
\item simplicial: $G_v\actson Y_v$ is a simplicial action on a simplicial tree;
\item of Seifert type: the vertex action $G_v\actson Y_v$ has kernel $N_v$, and the faithful action $G_v/N_v\actson Y_v$
is dual to an arational measured foliation on a closed $2$-orbifold with boundary;
\item axial: $Y_v$ is a line, and the image of $G_v$ in $\Isom(Y_v)$ is a finitely generated group 
acting with dense orbits on $Y_v$.
\end{enumerate*}
\end{thm*}

An \emph{infinite tripod} is the union of three semi-lines having a common origin $O$, and whose pairwise intersection is reduced to $\{O\}$.

If one assumes triviality of stabilizers of tripods and of unstable arcs (super-stability) in Main Theorem, 
one gets Sela's result.

The non-simplicial building blocks are canonical. Indeed, non-simplicial building blocks
are \emph{indecomposable} which implies that they cannot be split further into a graph of actions
(see Definition \ref{dfn_indecomposability} and Lemma \ref{lem_indec_component}).
If $T$ is not a line, simplicial building blocks can also be made canonical by imposing that each simplicial building block 
is an arc which intersects the set of branch points of $T$ exactly at its endpoints.

For simplicity, we did not state the optimal statement of Main Theorem, see Theorem \ref{thm_main} for more details.
In particular, one can say a little more about the tripod stabilizer on which $G$ splits.
This statement also includes a relative version: 
$G$ is only assumed to be finitely generated relative to a
finite set of elliptic subgroups $\calh=\{H_1,\dots,H_p\}$, and  
each $H_i$ is conjugate into a vertex group in the splittings of $G$ produced.\\

When $T$ has a decomposition into a graph of actions as in the conclusion of Main Theorem,
then $G$ splits over an extension of an arc stabilizer by a cyclic group except maybe if $T$ is a line. Thus we get:

\begin{UtiliseCompteurs}{corBF}
\begin{cor}
  Under the hypotheses of Main Theorem,
either $T$ is a line, or $G$ splits over a subgroup $H$ which is an extension
of a cyclic group by an arc stabilizer.
\end{cor}  
\end{UtiliseCompteurs}

\begin{rem*}
  In the conclusion of the corollary, $H$ is an extension by a \emph{full} arc stabilizer, 
and not of a subgroup fixing an arc. This contrasts with \cite[Theorem 9.5]{BF_stable}.
\end{rem*}

Here is a simple setting where the main theorem applies.
A group $H$ is \emph{small} (resp. \emph{slender}) if it contains no non-abelian free group
(resp. if all its subgroups are finitely generated).
The following corollary applies to the case where $G$ is hyperbolic relative to slender groups.

\begin{UtiliseCompteurs}{cor_petit}
\begin{cor}
  Let $G$ be a finitely generated group for which any small subgroup is finitely generated.
Assume that $G$ acts on an $\bbR$-tree $T$ with small arc stabilizers. 

Then 
either $G$ splits over the stabilizer of an unstable arc or over a tripod stabilizer,
or $T$ has a decomposition into a graph of actions as in Main Theorem.
In particular, $G$ splits over a small subgroup.
\end{cor}
\end{UtiliseCompteurs}

In some situations, one can control arc stabilizers in terms of tripod stabilizers. For instance, we get:

\begin{UtiliseCompteurs}{cor_tripod}
\begin{cor}
Consider a finitely generated group $G$ acting by isometries on an $\bbR$-tree $T$.
Assume that
\begin{enumerate*}
\item arc stabilizers have a nilpotent subgroup of bounded index (maybe not finitely generated);
\item tripod stabilizers are finitely generated  (and virtually nilpotent);
\item no group fixing a tripod is a proper subgroup of any conjugate of itself;
\item any chain $H_1\subset H_2\dots$
of subgroups fixing tripods stabilizes. 
\end{enumerate*}

Then  either $G$ splits over a subgroup having a finite index subgroup fixing a tripod, 
or $T$ has a decomposition as in the conclusion of Main Theorem.\\
\end{cor}
\end{UtiliseCompteurs}

Our proof relies on a particular case of Sela's Theorem, assuming triviality of arc stabilizers.
For completeness, we give a proof of this result in appendix \ref{sec_Sela}. 
To reduce to this theorem, we prove that under the hypotheses of Main Theorem, the action is \emph{piecewise stable},
meaning that any segment of $T$ is covered by finitely many stable arcs (see Definition \ref{dfn_stable}).
Piecewise stability implies that $T$ splits into a graph of actions where each vertex action
$G_v\actson T_v$ has trivial arc stabilizers up to some kernel (Theorem \ref{thm_pw2triv}).
Because $G_v$ itself need not be finitely generated, we need to extend Sela's result to 
finitely generated pairs (Proposition \ref{prop_sela_rel}).

So the main step in the proof consists in proving piecewise stability (Theorem \ref{thm_acc2pw}).
All studies of actions on $\bbR$-trees use resolutions by foliated $2$-complexes.
Such foliated $2$-complexes have a dynamical decomposition into two kinds of pieces:
\emph{simplicial} pieces where each leaf is finite, and \emph{minimal components} where every leaf is dense.
These resolutions give sequence of actions $G_k\actson T_k$ converging strongly to $T$, where $G$ is
the inductive limit of $G_k$.

Maybe surprisingly, the main difficulty in the proof arises from the simplicial pieces.
This is because minimal components
give large stable subtrees (Lemma \ref{lem_indec_stab}).
In particular, a crucial argument is a result saying that if all the groups $G_k$ split in some nice compatible way, then so does $G$.
In the case of free splittings, this is Scott's Lemma.
In our setting, we need an extended version of Scott's Lemma, which is of independent interest.

\begin{UtiliseCompteurs}{scott}
\begin{thm}[Extended Scott's Lemma]
Let $G_k\actson S_k$ be a sequence of non-trivial actions of finitely generated groups on simplicial trees,
and $(\phi_k,f_k):G_k\actson S_k\ra G_{k+1}\actson S_{k+1}$ be epimorphisms.
Assume that $(\phi_k,f_k)$ does not increase edge stabilizers in the following sense:
$$\forall e\in E(S_k),\forall e'\in E(S_{k+1}),\quad e'\subset f_k(e)\Rightarrow G_{k+1}(e')=\phi_k(G_k(e))$$

Then the inductive limit $G=\displaystyle\lim_{\ra} G_k$ has a non-trivial splitting over 
the image of an edge stabilizer of some $S_k$.
\end{thm}
\end{UtiliseCompteurs}

The paper is organized as follows.
Section \ref{sec_prelim} is devoted to preliminaries, more or less well known, 
except for the notion of indecomposability which seems to be new.
Section \ref{sec_scott} deals with Extended Scott's Lemma.
Section \ref{sec_acc2pw} proves piecewise stability.
Section \ref{sec_pw2triv} shows how piecewise stability and Sela's Theorem for actions with trivial arc stabilizers
allow to conclude.
Section \ref{sec_proof} gives the proof of the corollaries.
Section \ref{sec_example} contains our counterexample to Sela's misstated result.
Appendix \ref{sec_Sela} gives a proof of the version of Sela's result we need.





\section{Preliminaries}\label{sec_prelim}

\subsection{Basic vocabulary}

An $\bbR$-tree is a $0$-hyperbolic geodesic space.
In all this section, we fix an isometric action of a  group $G$ on an $\bbR$-tree $T$.
A subgroup $H\subset G$ is \emph{elliptic} if it fixes a point in $T$.
The action $G\actson T$ is \emph{trivial} if $G$ is elliptic.
The action $G\actson T$ is \emph{minimal} if $T$ has no proper $G$-invariant subtree.
We consider the trivial action of $G$ on a point as minimal.
When $G$ contains a hyperbolic element, then $T$ contains a unique minimal $G$-invariant subtree
$\Tmin$, and $\Tmin$ is the union of  axes of all hyperbolic elements.

An \emph{arc} is a set homeomorphic to $[0,1]$.
A subtree is \emph{non-degenerate} if it contains an arc.

Say that an action on an $\bbR$-tree $G\actson T$ is \emph{simplicial} 
it can be obtained from an action on a combinatorial tree
by assigning equivariantly a positive length for each edge.
If $S$ is a simplicial tree, we denote by $V(S)$ its set of vertices, and $E(S)$ its set of oriented edges.

A \emph{morphism of $\bbR$-trees} $f:T\ra T'$ is a $1$-Lipschitz map
such that any arc of $T$ can be subdivided into a finite number of sub-arcs
on which $f$ is isometric.

For basic facts about $\bbR$-trees, see \cite{Sh_dendrology,Chi_book}.

\subsection{Stabilities}
Stability hypotheses say how arc stabilizers are nested. 
 By stabilizer, we always mean \emph{pointwise} stabilizer. 
When we will talk about the global stabilizer, we will mention
it explicitly.
If $G$ acts on $T$, and $X\subset T$, we will denote by $G(X)$ the (pointwise)
stabilizer of $X$.

\begin{dfn}[Stabilities]\label{dfn_stable} 
Consider an action of a group $G$ on an $\bbR$-tree $T$.
\begin{itemize*}
\item A non-degenerate subtree $Y\subset T$ is called \emph{stable} if
  for every arc $I\subset Y$, $G(I)=G(Y)$.
  
\item $T$ is \emph{BF-stable} (in the sense of Bestvina-Feighn \cite{BF_stable})
  if every   arc of $T$ contains a stable arc. 

\item $T$ satisfies \emph{the ascending chain condition} if for any
sequence of arcs $I_1\supset I_2\supset \dots$ whose lengths converge to $0$,
the sequence of stabilizers $G(I_1)\subset G(I_2)\subset\dots$ is eventually constant.

\item $T$ is \emph{piecewise stable} if any
  arc of $T$ can covered by finitely many stable arcs. Equivalently,
  the action is piecewise stable if for all $a\neq b$, the arc $[a,b]$
  contains a stable arc of the form $[a,c]$ (one common endpoint).
  
\item $T$ is \emph{super-stable} if any  arc
 with non-trivial stabilizer is stable.
\end{itemize*}
\end{dfn}

\begin{rem*}
  The first definition of piecewise stability clearly implies the second one. 
Conversely, assume that any arc $[a,b]$ contains a stable arc of the form $[a,c]$.
Let $I$ be any arc of $T$. Thus any point of $I$ has a neighbourhood in $I$ which is the union of
at most two stable arcs. By compactness, $I$ is covered by finitely many stable arcs.
\end{rem*}

\begin{rem*}
  Clearly, super-stability or piecewise stability imply the ascending chain condition
which implies BF-stability.
If $T$ is a stable subtree of itself, then any arc stabilizer fixes $T$; in other words,
if $N$ is the kernel of the action $G\actson T$ (\ie $N$ is the set of elements acting as the identity),
$G/N$ acts with trivial arc stabilizers on $T$.
\end{rem*}

\subsection{Graphs of actions on $\bbR$-trees and transverse coverings}
\label{sec_goa}

Graphs of actions on $\bbR$-trees are a way of gluing equivariantly $\bbR$-trees
(see \cite{Skora_combination} or \cite{Lev_graphs}).
Here, we rather follow \cite[section 4]{Gui_limit}.

\begin{dfn}
  A \emph{graph of actions on $\bbR$-trees} $\calg=(S,(Y_v)_{v\in V(S)},(p_e)_{e\in E(S)})$  consists of the
  following data:
\begin{itemize*}
\item a simplicial tree $S$ called the \emph{skeleton} with a
  simplicial action without inversion of a group $G$;
\item for each vertex $v\in S$, an $\bbR$-tree $Y_v$ (called vertex
  tree or vertex action); 
\item for each oriented edge $e$ of $S$ with terminal vertex v, an
  \emph{attaching point} $p_e\in Y_v$.
\end{itemize*}
All this data should be invariant under $G$: $G$ acts on the
disjoint union of the vertex trees so that the projection 
$Y_v\mapsto v$ is equivariant;
and for every $g\in G$, $p_{g.e}=g.p_e$. 
\end{dfn}

\begin{rem*}
Some vertex trees may be reduced to a point (but they are not allowed to be empty).
  The definition implies that $Y_v$ is $G(v)$-invariant
and that $p_e$ is $G(e)$-invariant. 
\end{rem*}

To a graph of actions $\calg$ corresponds an 
$\bbR$-tree $T_\calg$ with a natural action of $G$. 
Informally, $T_\calg$ is obtained from the disjoint union of the vertex trees by
identifying the two attaching points of each edge of $S$ (see \cite{Gui_limit} for a formal definition). 

Alternatively, one can define a graph of actions as a graph of groups $\Gamma$
with an isomorphism $G\simeq\pi_1(\Gamma)$
together with an action of each vertex group $G_v$ on an $\bbR$-tree $Y_v$ and
for each oriented edge $e$ with terminus vertex $v$, a point of $Y_v$ fixed by the image of $G_e$ in $G_v$.
This is where the terminology comes from.

\begin{dfn}
Say that $T$ splits as a graph of actions $\calg$ if
there is an equivariant isometry between 
  $T$ and  $T_\calg$.
\end{dfn}

Transverse coverings are very convenient when working with graphs of actions. 

\begin{dfn}
A  \emph{transverse covering} of an $\bbR$-tree $T$ is a covering of $T$
by a family of subtrees $\caly=(Y_v)_{v\in V}$ such that
\begin{itemize*}
\item every $Y_v$ is a closed subtree of $T$
\item every arc of $T$ is covered by finitely many subtrees of $\caly$
\item for $v_1\neq v_2\in V$, $Y_{v_1}\cap Y_{v_2}$ contains at most one point
\end{itemize*}
\end{dfn}

When $T$ has an action of a group $G$, we always require the family $\caly$ to be $G$-invariant.

\begin{rem*}
In the definition above, if some subtrees $Y_v$ are reduced to a point, we may as well
forget them in $\caly$.
Moreover, given a covering by subtrees $Y_v$ which are not closed, but
  which satisfy the two other conditions for a transverse covering, then the
  family $(\ol{Y_v})_{v\in V}$ of their closure is a transverse covering.
\end{rem*}

The relation between graphs of actions and transverse coverings is contained in the following result.

\begin{lem}[{\cite[Lemma 4.7]{Gui_limit}}]\label{lem_transverse_cov}
Assume that $T$ splits as a graph of actions with vertex trees $(Y_v)_{v\in V(S)}$.
Then the family $\caly=(Y_v)_{v\in V(S)}$ is a transverse covering of $T$.

Conversely, if $T$ has a transverse covering by a family $\caly=(Y_v)_{v\in V}$ of non-degenerate trees,
then $T$ splits as a graph of actions whose non-degenerate vertex trees are the $Y_v$.
\end{lem}

We recall for future use the definition of the graph of actions induced by a transverse covering $\caly$.
We first define its skeleton $S$.
Its vertex set $V(S)$ is $V_0(S)\cup V_1(S)$ where $V_1(S)$ is the set of subtrees $Y\in\caly$, 
and $V_0(S)$ is the set of points $x\in T$ lying 
in the intersection of two distinct subtrees in $\caly$.
There is an edge $e=(x,Y)$ between $x\in V_0(S)$ and $Y\in V_1(S)$ if $x\in Y$. 
The vertex tree of $v\in V(S)$ is the corresponding subtree of $T$ (reduced to a point if and only if $v\in V_0(S)$).
The two attaching points the edge $(x,Y)$ are the copies of $x$ in $\{x\}$ and $Y$ respectively.

\begin{example*}
  Let $G\actson T$ be a simplicial action. Then the family of its edges is a transverse covering of $T$.
If $T$ has no terminal vertex, then 
the skeleton of this transverse covering is the barycentric subdivision of $T$.
\end{example*}

Here is another general example which will be useful.

\begin{lem}\label{lem_decompo_minimal}
  Assume that $G\actson T$ contains a hyperbolic element and denote by
 $\Tmin$ be the minimal subtree of $T$.
Let $\caly_0$ be the set of closures of connected components of $T\setminus \ol\Tmin$.

Then $\{\ol\Tmin\}\cup \caly_0$ is a transverse covering of $T$.
\end{lem}

\begin{rem*}
If $Y_v\in\caly_0$, then its global stabilizer $G_v$ fixes the point
$Y_v\cap \ol\Tmin$. In particular, this lemma says that a finitely supported action (in the sense of Definition \ref{dfn_span})
can be decomposed into a graph of actions where each vertex action
is either trivial or has a dense minimal subtree.
\end{rem*}

\begin{proof}
Any arc $I\subset T$ is covered by $\ol\Tmin$ and at most two elements of $\caly_0$, containing the endpoints of $I$.
The other properties of the transverse covering are clear.
\end{proof}

Consider a graph of actions $\calg=(S,(Y_v)_{v\in V(S)},(p_e)_{e\in E(S)})$, and $T_\calg$ the corresponding $\bbR$-tree.
Given a $G$-invariant subset $V'\subset V(S)$, we want to define an $\bbR$-tree by collapsing all the trees $(Y_v)_{v\in V'}$.
For each $v\in V'$, replace $Y_v$ by a point, and for each edge $e$ incident on $v\in V'$, change $p_e$ accordingly.
Let $\calg'=(S,(Y'_v)_{v\in V(S)},(p'_e)_{e\in E(S)})$ the corresponding graph of actions, and $T_\calg'$ the corresponding $\bbR$-tree.

\begin{dfn}\label{dfn_collapse}
We say that the tree $T_{\calg'}$ is obtained from $T_\calg$ by collapsing the trees $(Y_v)_{v\in V'}$.
\end{dfn}

Let $p_v:Y_v\ra Y'_v$ the natural map (either the identity or the constant map),
and $p:T_\calg\ra T_\calg'$ the induced map.

\begin{dfn}
  We say that a map $f:T\ra T'$ \emph{preserves alignment} if the three
  following equivalent conditions hold:
  \begin{itemize}
  \item for all $x\in [y,z]$, $f(x)\in [f(y),f(z)]$.
  \item the preimage of a convex set is convex,
  \item the preimage of a point is convex.
  \end{itemize}
\end{dfn}
See for instance \cite[Lemma 1.1]{Gui_coeur} for the equivalence.

\begin{lem}\label{lem_collapse}
  The map $p:T_\calg\ra T_\calg'$ preserves alignment.
In particular, if $T_\calg$ is minimal, so is $T_\calg'$.
\end{lem}

\begin{proof}
Consider a point $x\in T_{\calg'}$. 
Let $E'_v=\{x\}\cap Y'_v$.
If $u,v$ are such that $E'_{u},E'_v\neq\es$ for $i=1,2$,
there is a path $u=v_0,\dots,v_n=v$  such
that the attaching points of the edge $v_iv_{i+1}$ coincide with the point in $E'_{v_i}$
and $E'_{v_{i+1}}$ (in particular $E'_{v_i}\neq \es$).

Consider $E_v=p_v\m(E'_v)$.
For each $v$, either $E_v=Y_v$ or $E_v$ consists of at most one point.
The existence of the path above shows that the image of $E_v$ in $T_\calg$ is connected.
This proves that $p\m(\{x\})$ is convex so $p$ preserves alignment.

Let $Y'\subset T_\calg'$ be a non-empty $G$-invariant subtree.
Then $p\m(Y')$ is a non-empty  $G$-invariant subtree, so $p\m(Y')=T_\calg$.
Since $p$ is onto, $Y'=T_{\calg'}$, and minimality follows.
\end{proof}



\subsection{Actions of pairs, finitely generated pairs}

By a pair of groups $(G,\calh)$, we mean a group $G$ together with a finite family of subgroups
$\calh=\{H_1,\dots,H_p\}$. We also say that the groups $H_i$ are peripheral subgroups of $G$.
In fact, each peripheral subgroup is usually only defined up to conjugacy so it would be
more correct to define $\calh$ as a finite set of conjugacy class of subgroups.

An action of a pair $(G,\calh)$ on a tree is an action of $G$ in which each subgroup $H_i$ is elliptic.
When the tree is simplicial, we also say that this is a splitting of $(G,\calh)$ or a splitting 
of $G$ \emph{relative} to $\calh$.

\begin{dfn} \label{dfn_relfg}
  A pair $(G,\calh)$ is \emph{finitely generated} if there exists a finite set $F\subset G$ such that 
$F\cup H_1\cup\dots\cup H_p$
generates $G$.
We also say in this case that $G$ is finitely generated relative to $\calh$.
\end{dfn}

We say that $F$ is a generating set of the pair $(G,\calh)$.

\begin{rem*}
If $(G,\{H_1,\dots,H_p\})$ is finitely generated, then for any $g_i\in G$, so is
$(G,\{H_1^{g_1},\dots,H_p^{g_p}\})$: just add $g_i$ to $F$.
Of course, if $G$ is finitely generated, then $(G,\calh)$ is finitely generated for any $\calh$.
\end{rem*}

By Serre's Lemma, if $(G,\calh)\actson T$ is a non-trivial action of a finitely generated pair, then
$G$ contains a hyperbolic element (\cite[Proposition 6.5.2]{Serre_arbres}).

If $G_v$ is a vertex group in a finite graph of groups $\Gamma$,
then it has a natural structure of pair $(G_v,\calh_v)$ consisting of (representatives of)
conjugacy classes of the images in $G_v$ of the incident edge groups. 
We call this pair the \emph{peripheral structure} of $G_v$ in $\Gamma$.

\begin{lem} \label{lem_relfg}
Consider a graph of groups $\Gamma$ and $G=\pi_1(\Gamma)$.

If $\Gamma$ is finite and $G$ is finitely generated, 
then the peripheral structure $(G_v,\calh_v)$ of each vertex group $G_v$ is finitely generated.
\end{lem}

\begin{proof}
  We use the notations of \cite{Serre_arbres}.  Choose a maximal
  subtree $\tau$ of $\Gamma$. For each vertex $v$, the vertex
  group $G_v$ is now identified with a subgroup of
  $\pi_1(\Gamma,\tau)$. Let $F$ be a finite generating set of $\pi_1(\Gamma,\tau)$.
  Since we can write each $g\in F$ as a product of edges of $\Gamma$ and of elements of the
  vertex groups, there is a generating set of $\pi_1(\Gamma,\tau)$
  consisting of the edges of $\Gamma$ and of a finite set $F'$ of elements of the vertex groups of
  $\Gamma$.

For each vertex $v\in\Gamma$, let $H_v\subset G_v$ be the subgroup
generated by the elements of $F'$ which lie in $G_v$, and by the image in $G_v$ of incident edge groups. Of
course, $\pi_1(\Gamma,\tau)$ is generated by the groups $H_v$ and the
edges of $\Gamma$.

Fix $v\in\Gamma$. We shall prove that $G_v=H_v$.
Take $g\in G_v$,  and write $g$ as $g=g_0e_1g_1e_2...e_kg_k$ where $v_0 e_1
v_1\dots e_k v_k$ is a loop based at $v$, and each $g_i$ is in
$H_{v_i}$.  If the word has length $0$, there is nothing to prove.  If
not, we shall find a shorter word of this form
representing $g$, and the lemma will be proved by induction. 
If this word has not length $0$, it is not 
reduced as a word in the graph of groups.
Therefore, there are two
consecutive edges $e_i,e_{i+1}$ such that $e_{i+1}=\Bar{e_i}$, and
$g_i\in j_{e_i}(G_{e_i})$ (where $j_e$ denotes the edge morphisms
of $\Gamma$).  Thus, we can write $h=g_0e_1\dots
e_{i-1}g_{i-1}g'_ig_{i+1}e_{i+2}\dots e_kg_k$ where $g'_i\in
j_{e_{i+1}}(G_{e_{i+1}})$.  Note that $g_{i-1}g'_ig_{i+1}\in
H_{v_{i-1}}=H_{v_{i+1}}$ hence we found a shorter word of the required
form representing $g$.
\end{proof}

\begin{lem}\label{lem_relfg_rel}
  Let $(G,\calh)$ be a finitely generated pair, and consider a (relative) splitting of $(G,\calh)$.
Let $v$ be a vertex of the corresponding graph of groups.

Then $G_v$ is finitely generated relative to a family $\calh_v$ consisting of 
the images in $G_v$ of the incident edge groups together with at most one conjugate of each $H_i$.
\end{lem}

\begin{proof}
Let $G\simeq\pi_1(\Gamma)$ be a splitting of $G$ as a graph of groups.
We can assume that $H_i\subset G_{v_i}$ for some vertex $v_i$ of $\Gamma$.
For each index $i$, consider a finitely generated group $\Hat H_i$ containing $H_i$.
Consider the graph of groups $\Hat \Gamma$ obtained from $\Gamma$ 
by adding for each $i$ a new edge $e_i$ edge carrying $H_i$ incident on  a new vertex $\Hat v_i$ carrying $\Hat H_i$,
and by gluing the other side of $e_i$ on $v_i$.
Thus, $e_i$ carries the amalgam $G_{v_i}*_{H_i} \Hat H_i$.
The fundamental group of  $\Hat\Gamma$ is finitely generated.
By Lemma \ref{lem_relfg}, the peripheral structure of $G_v$ in $\Hat G_v$ is a finitely generated pair.
The lemma follows.
\end{proof}

\subsection{Finitely supported actions}

A \emph{finite tree} in an $\bbR$-tree is the convex hull of a finite set.

\begin{dfn}\label{dfn_span}
A set $K\subset T$ \emph{spans} $T$ if
any arc $I\subset T$ is covered by finitely many translates of $K$.
The tree $T$ is \emph{finitely supported} if it is spanned by a finite tree $K$.
\end{dfn}

If $T$ is a simplicial tree, then $T$ is finitely supported if and only if $T/G$ is a finite graph.
Any minimal action of a finitely generated group is finitely supported.
More generally we have:

\begin{lem}
A minimal action of a finitely generated pair on an $\bbR$-tree is finitely supported.
\end{lem}

\begin{proof}
Consider $(G,\{H_1,\dots,H_p\})$ a finitely generated pair acting minimally on $T$.
Consider a finite generating set $\{f_1\dots,f_q\}$ of the pair $(G,\{H_1,\dots,H_p\})$.
Let $x_i\in T$ be a point fixed by $H_i$ and $x\in T$ be any point.
Let $K$ be the convex hull of $\{x,f_1.x,\dots,f_q.x,x_1,\dots,x_p\}$.
Then $G.K$ is connected because for each generator $g\in \{f_1,\dots,f_q\}\cup H_1\cup\dots \cup H_p$,
$g.K\cap K\neq\es$. By minimality, $G.K=T$ and $T$ is finitely supported.
\end{proof}

\begin{lem}\label{lem_supp_fini}
Consider an action of a finitely generated pair $(G,\calh)\actson T$,
$\caly$ a transverse covering of $T$, and $S$  the skeleton of $\caly$.

Then any $H\in\calh$ is elliptic in $S$. 
Moreover, if $G\actson T$ is minimal (resp. finitely supported)
then so is $G\actson S$.
\end{lem}

\begin{proof}
The statement about minimality is proved in \cite[Lemma 4.9]{Gui_limit}.

We use the description of the skeleton given after Lemma \ref{lem_transverse_cov}.
We can assume that every subtree $Y\in\caly$ is non-degenerate.
Let $K$ be a finite tree spanning $T$. It is covered by finitely many trees of $\caly$.
It follows that $\caly/G=V_1(S)/G$ is finite.
Choose a lift of each element of $V_1(S)/G$ in $V_1(S)$,
and consider the convex hull $L\subset S$ of those vertices.
Let $x_i\in T$ be a point fixed by $H_i$.
Each $H_i\in \calh$ fixes a vertex $v_i\in S$, lying in $V_0(\caly)$ or $V_1(\caly)$
according to whether $x_i$ lies in two distinct elements of $\caly$ or not.

Let $\{f_1,\dots,f_q\}$ be a finite generating set of $(G,\calh)$.
Let $L_0$ be the convex hull of $L\cup f_1.L\cup\dots \cup f_q.L\cup\{v_1,\dots,v_p\}$. 
The set $S_0=G.L_0$ is a subtree of $S$ containing $V_1(S)$. It is finitely supported.
If $S_0\neq S$, consider $x\in V_0(S)\setminus S_0$. 
Recall that the vertex $x\in V_0(S)$ corresponds to a point $x\in T$.
Any edge of $S$ incident on $x$ corresponds to a subtree $Y\in V_1(S)$ such that $x\in Y$.
Since $S_0$ is a subtree and contains $V_1(S)$, there is exactly one such edge.
If follows that $x$ belongs to exactly one $Y\in\caly$, contradicting the definition of $V_0(S)$.
Therefore, $S=S_0$ is finitely supported.
\end{proof}

If $T$ splits as a graph of actions, 
vertex actions may fail to be minimal, even if $T$ is minimal.
However, the property of being finitely supported is inherited by vertex actions:

\begin{lem}\label{lem_support}
Assume that $T$ splits as a graph of actions and that $T$ is finitely supported.
Then each vertex action is finitely supported.
\end{lem}

\begin{proof}
  Let $K$ be a finite tree spanning of $T$.
Let $\caly$ be a transverse covering of $T$, and $Y_0\in\caly$ be a vertex tree.
Consider a finite subset $F$ of $\caly$ which covers $K$,
and let $\{g_1.Y_0,\dots,g_p.Y_0\}$ be the set of elements of $F$ lying in the orbit of $Y_0$.
Consider the convex hull $K_0$  of the finite trees $Y_0\cap g_i\m K$ for $i\in\{1,\dots, p\}$.
It is easy to check that the finite tree $K_0$ spans $Y_0$.
\end{proof}





\subsection{Indecomposability}

Indecomposability is a slight modification of the \emph{mixing} property
introduced by Morgan in \cite{Mo_ergodic}.

\begin{figure}[htbp]
  \begin{center}
    \includegraphics{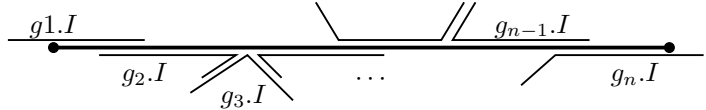}
\caption{Indecomposability\label{fig_indecomposability}.}
  \end{center}
\end{figure}

\begin{dfn}[Indecomposability] \label{dfn_indecomposability}
 A non-degenerate subtree  $Y\subset T$ is called \emph{indecomposable} if for every pair of arcs
  $I,J\subset Y$, there is a finite sequence
  $g_1.I,\dots,g_n.I$ which covers $J$ and such that
  $g_i.I\cap g_{i+1}.I$ is non-degenerate (see figure
  \ref{fig_indecomposability}).
\end{dfn}

\begin{rem*}
In the definition of indecomposability,
one cannot assume in general that $g_i.I\cap J$ is non-degenerate for all $i$.
This is indeed the case if $Y=T$ is dual to a measured foliation on a surface 
having a 4-pronged singularity, 
$J$ is an arc in $T$ represented by a small transverse segment containing this singularity
and joining two opposite sectors of the singularity, and $I$ is represented 
by a transverse segment disjoint from the singularities.
\end{rem*}

The following property explains the choice of the terminology.
\begin{lem} \label{lem_indec_component}
  If $Y\subset T$ is indecomposable, and if $T$ splits as a graph of
  actions, then $Y$ is contained in a vertex tree.
\end{lem}

\begin{proof}
  Consider an arc $I$ contained in $Y\cap Y_v$ for some vertex tree $Y_v$ of the
  decomposition.
  Consider an arc $J\subset Y$, and a finite sequence
  $g_1.I,\dots,g_n.I$ which covers $J$ and such that
  $g_i.I\cap g_{i+1}.I$ is non-degenerate.  We have
  $g_1.I\subset g_1.Y_v$, but $g_1.Y_v\cap
  g_2.Y_v$ is non-degenerate so $g_1.Y_v=g_2.Y_v$. By
  induction, we get that all the translates $g_i.I$ lie in the
  same vertex tree $g_1.Y_v$, so that $J$ lies in a vertex tree.
  Since this is true for every arc $J\subset Y$, the lemma follows.
\end{proof}

\begin{lem}\label{lem_indecompo} 
  \begin{enumerate*}
  \item If $f:T\ra T'$ is a morphism of $\bbR$-trees, and if $Y\subset
    T$ is indecomposable, then so is $f(Y)$.
  \item If $Y\subset T$ is indecomposable, then the orbit of any point
    $x\in Y$ meets every arc $I\subset Y$ in a dense subset.
  \item If $T$ itself is indecomposable, then it is minimal (it has no non-trivial invariant subtree).
  \item Assume that $(Y_v)_{v\in V}$ is a transverse covering of $T$. Let $Y_{v_0}$ be a vertex tree and $H$ its global stabilizer.
If $Y_{v_0}$ is an indecomposable subtree of $T$, then $Y_{v_0}$ is indecomposable for the action of $H$.
  \end{enumerate*}
\end{lem}

\begin{proof}
Let $I=[f(a),f(b)]$ and $J=[f(x),f(y)]$ in $f(Y)$.
Choose $I'\subset [a,b]$ so that $f(I')\subset I$, and
let $J'=[x,y]$, so $f(J')$ contains $J$.
Indecomposability of $Y$ now clearly implies indecomposability of $f(Y)$.


Statement 2 follows from the fact that
for any arc $I'\subset I$, there exists $g_1$ such that 
$x\in g_1.I'$, so $g_1\m.x\in I'$.

Let's prove statement 3.
It follows from statement 2 that $G$ contains a hyperbolic element.
Consider $\Tmin$ the minimal $G$-invariant subtree.
By statement 2, every orbit meets $\Tmin$; if follows that $T=\Tmin$.

For statement 4, consider $I,J\subset Y_{v_0}$, and 
$g_1,\dots,g_n\in G$ such that $J\subset g_1.I\cup \dots \cup g_n.I$ with $g_i.I\cap g_{i+1}.I$ non-degenerate.
In particular, $g_i.Y_{v_0}\cap g_{i+1}.Y_{v_0}$ is non-degenerate, so $g_i.Y_{v_0}= g_{i+1}.Y_{v_0}$,
and $g_{i+1}g_i\m\in H$.
Consider an index $i$ such that  $g_i.I\cap J$ is non-degenerate.
Since $J\subset Y_{v_0}$,  $g_i.Y_{v_0}\cap Y_{v_0}$ is non-degenerate, so $g_i\in H$.
Therefore, $g_1,\dots, g_n\in H$.

\end{proof}

The conjunction of indecomposability and stability has a nice
consequence:
\begin{lem}\label{lem_indec_stab}
  If $T$ is BF-stable, and if $Y\subset T$ is
  indecomposable, then $Y$ is a stable subtree.
\end{lem}

\begin{proof}
First, one easily checks that if $K_1,K_2$ are two stable subtrees of $T$ such that
$K_1\cap K_2$ is non-degenerate, then $K_1\cup K_2$ is a stable subtree.
Consider a stable arc $I\subset Y$, and any other arc $J\subset Y$.
Since $J$ is contained in a union $g_1.I\cup\dots\cup g_n.I$ with $g_i.I\cap g_{i+1}.I$ non-degenerate,
 $J$ is stable.
Since this holds for any arc $J$, $Y$ is a stable subtree of $T$.
\end{proof}

\subsection{Geometric actions and strong convergence}

We review some material from \cite{LP} where more details can be found.

\subsubsection*{Strong convergence}

Given two actions $G\actson T$, $G'\actson T'$,
we write $(\phi,f):G\actson T\ra G'\actson T'$
when $\phi:G\ra G'$ is a morphism,
and $f:T\ra T'$ is a $\phi$-equivariant map.
We say that $(\phi,f)$ is onto if both $f$ and $\phi$ are.
We say that $(\phi,f)$ is morphism of $\bbR$-trees if $f$ is.

A \emph{direct system} of actions on $\bbR$-trees
is a sequence of actions of finitely generated groups $G_k\actson T_k$ and an action $G\actson T$,
with surjective morphisms of $\bbR$-trees $(\phi_k,f_k):G_k\actson T_k\onto G_{k+1}\actson T_{k+1}$
and  $(\Phi_k,F_k):G_k\actson T_k\onto G\actson T$ such that the
  following diagram commutes:
  $$\xymatrix@1@R=0.5cm{T_k\ar[r]_{f_k} \ar@/^0.9cm/[rrr]|{F_k}
    \ar@(dl,dr)[]
    & T_{k+1}  \ar@/^0.5cm/[rr]|{F_{k+1}} \ar@(dl,dr)[] 
    & \cdots
    & T \ar@(dl,dr)[] \\
    G_k\ar[r]^{\phi_k} \ar@/_0.7cm/[rrr]|{\Phi_k}
    &G_{k+1} \ar@/_0.4cm/[rr]|{\Phi_{k+1}} 
    &\cdots & G}
  $$

For convenience, we will use the notation $f_{kk'}=f_{k'-1}\circ\dots\circ f_k:T_k\ra T_{k'}$
and $\phi_{kk'}=\phi_{k'-1}\circ\dots\circ \phi_k:G_k\ra G_{k'}$.

  \begin{dfn}[Strong convergence]
  A direct system of finitely supported actions of
  groups on $\bbR$-trees $G_k\actson T_k$
  \emph{converges strongly} to $G\actson T$ if
\begin{itemize*}
\item $G$ is the direct limit of the groups $G_k$
\item for every finite tree $K\subset T_k$, there exists $k'\geq k$ such
  that $F_{k'}$ restricts to an isometry on $f_{kk'}(K)$,
\end{itemize*}
The action  $G\actson T$ is the \emph{strong limit} of this direct system.
  \end{dfn}

In practice, strong convergence allows to lift the action of a finite number of elements on a finite subtree of $T$
to $T_k$ for $k$ large enough. 

As an application  of the definition, we prove the following useful lemma.
\begin{lem}\label{lem_strong_minimal}
  Assume that $G_k\actson T_k$  converges strongly to $G\actson T$.
Assume that $G\actson T$ is minimal.

Then for $k$ large enough, $G_k\actson T_k$ is minimal.
\end{lem}

\begin{proof}
By definition of strong convergence, $T_0$ is finitely supported.
  Let $K_0\subset T_0$ be a finite tree spanning $T_0$.
Let $K$ (resp. $K_k$) be the image of $K_0$ in $T$ (resp. in $T_k$).
Choose some hyperbolic elements $g_1,\dots,g_p\in G$ whose axes cover $K$.
Let $K'$ be the convex hull of $K\cup g_1.K\cup\dots\cup g_p.K$.
Choose $ g_1^0,\dots g_p^0\in G_0$ some preimages of $g_1,\dots g_p$
and let $g_i^k$ be the image of $g_i^0$ in $G_k$.
Let $K'_0\subset T_0$ be the convex hull  of $K_0 \cup g_1^0.K_0\cup\dots\cup g_p^0.K_0$.
Let $K'_k$ be the image of $K'_0$ in $T_k$.
Take $k$ large enough so that $F_k$ induces an isometry between $K'_k$ and $K'$.
Then the axes of $g_1^k,\dots,g_p^k$ cover $K_k$. Since $K_k$ spans $T_k$,
$T_k$ is a union of axes, so $T_k$ is minimal.
\end{proof}

\subsubsection*{Geometric actions}

A \emph{measured foliation} $\calf$ on a $2$-complex $X$
consists of the choice, for each closed simplex $\sigma$ of $X$ of a non-constant affine map
$f_\sigma:\sigma\ra \bbR$ defined up to post-composition by an isometry of $\bbR$, in a way which  
is consistent under restriction to a face: if $\tau$ is a face of $\sigma$, then 
$f_\tau=\phi\circ(f_\sigma)_{|\tau}$ for some isometry of $\phi$ of $\bbR$.
Level sets of $f_\sigma$ define a foliation on each closed simplex. 
Leaves of the foliations on $X$ are defined as the equivalence classes of the equivalence relation generated 
by the relation \emph{$x,y$ belong to a same closed simplex $\sigma$ and $f_\sigma(x)=f_\sigma(y)$}.

This also defines a transverse measure as follows:
the transverse measure $\mu(\gamma)$ of a path $\gamma:[0,1]\ra \sigma$ 
to the foliation is the length of the path $f_\sigma\circ\gamma$. 
The transverse measure of a path which is a finite concatenation of paths contained in simplices
is simply the sum of the transverse measures of the pieces.
The transverse measure is invariant under the holonomy along the leaves.
We also view this transverse measure as a metric on each transverse edge.

For simplicity, 
we will say \emph{foliated 2-complex} to mean a $2$-complex endowed with a measured foliation.

The pseudo-metric 
$$\delta(x,y)=\inf\{\mu(\gamma)\text{ for $\gamma$ joining $x$ to $y$}\}$$
is zero on each leaf of $X$. 
By definition, the \emph{leaf space made Hausdorff} $X/\calf$ of $X$ is the metric space obtained from $X$
by making $\delta$ Hausdorff, \ie by identifying points at pseudo-distance $0$. 
In nice situations, $\delta(x,y)=0$ if and only if $x$ and $y$ are on the same leaf.
In this case, $X/\calf$ coincides with the space of leaves of the foliation,
and we say that the \emph{leaf space is Hausdorff}.


\begin{thm}[{\cite[Proposition 1.7]{LP}}] \label{thm_LP_leaves}
Let $(X,\calf)$ be a foliated $2$-complex.
  Assume that $\pi_1(X)$ is generated by free homotopy classes of curves contained in leaves.

Then $X/\calf$ is an $\bbR$-tree.
\qed
\end{thm}


\begin{dfn}[Geometric action]\label{dfn_dual}
An action of a finitely generated group $G$ on an $\bbR$-tree $T$ 
is \emph{geometric} if there exists a foliated $2$-complex $X$ endowed with a free, properly discontinuous, cocompact action
of $G$, such that 
\begin{itemize*}
\item each transverse edge of $X$ isometrically embeds into $X/\calf$
\item $T$ and $X/\calf$ are equivariantly isometric.
\end{itemize*}
In this case, we say that $T$ is \emph{dual} to $X$.
\end{dfn}

\begin{rem*}
  The definition in \cite{LP} is in terms of a compact foliated $2$-complex $\Sigma$ (here $\Sigma=X/G$)
and of a Galois covering $\ol \Sigma$ with deck group $G$ (here $\ol\Sigma=X$).
The two points of view are clearly equivalent.

This definition requires $G$ to be finitely generated. For finitely generated pairs, one should 
weaken the assumption of cocompactness, but we won't enter into this kind of consideration.
\end{rem*}

\subsubsection*{Decomposition of geometric actions}

\begin{prop}\label{prop_decompo}
  Let $G\actson T$ be a geometric action dual to a $2$-complex $X$ whose fundamental group is generated
by free homotopy classes of curves contained in leaves. 

Then $T$ has a decomposition into a graph of actions where each non-degenerate vertex action is either 
indecomposable, or is an arc containing no branch point of $T$ except possibly at its endpoints.
\end{prop}

The arcs in the lemma above are called \emph{edges},
and the indecomposable vertex actions are called \emph{indecomposable components}.

The basis of the proof of this proposition is a version of a theorem by Imanishi
giving a dynamical decomposition of a compact foliated $2$-complex $\Sigma$ (\cite{Imanishi}).
A leaf of $\Sigma$ is \emph{regular} if it contains no vertex of $\Sigma$. 
More generally, a leaf segment (\ie a path contained in a leaf) is \emph{regular} if it contains no vertex of $\Sigma$.
A leaf or leaf segment which is not regular is \emph{singular}.

Let $\Sigma^*=\Sigma\setminus V(\Sigma)$ be the complement of the vertex set of $\Sigma$.
It is endowed with the restriction of the foliation of $\Sigma$.
Let $C^*\subset \Sigma^*$ be the
union of leaves of $\Sigma^*$ which are closed but not compact.

\begin{dfn}\label{dfn_cut}
We call the set $C=C^*\cup V(\Sigma)$ \emph{the cut locus} of $\Sigma$.  
\end{dfn}

The cut locus is a finite union of leaf segments joining two vertices of $\Sigma$.
In particular, $C$ is compact, and $\Sigma\setminus C$ consists of finitely many connected components.
Each component of  $\Sigma\setminus C$ is a union of leaves of $\Sigma^*$.

The dynamical decomposition of $\Sigma$ is as follows:
\begin{prop}[{\cite[Proposition 1.5]{LP}}]
  Let $U$ be a component of $\Sigma\setminus C$. Then either every leaf of $U$ is compact,
or every leaf of $U$ is dense in $U$.\qed
\end{prop}

\begin{proof}[Proof of Proposition \ref{prop_decompo}]
  Let $X$ be a foliated $2$-complex such that $T$ is dual to $X$.
Denote by $q:X\ra T$ the natural map (recall that $T$ is the leaf space made Hausdorff of $X$).
The quotient $\Sigma=X/G$ is a compact foliated $2$-complex
and the quotient map $\pi:X\ra \Sigma$ is a covering map.
Let $C$ be the cut locus of $\Sigma$ and $\Tilde C$ its preimage in $X$. 

Let $(U_v)_{v\in V}$ be the family of connected components of $X\setminus \Tilde C$.
Let $\ol U_v$ be the closure of $U_v$ in $X$, and $Y_v=q(\ol U_v)\subset T$.
Note that $\ol U_v\setminus U_v\subset \Tilde C$ is contained in a union of singular leaves
and that any leaf segment in $U_v$ is regular.
We shall first prove that the family $\caly=(Y_v)_{v\in V}$ is a transverse covering of $T$.
We will need the following result.

\begin{prop}[{\cite[Lemma 3.4]{LP}}]\label{prop_LP_separation}
Assume that $\pi_1(X)$ is generated by free homotopy classes of curves contained in leaves.

Then there exists a countable union of leaves $\cals$ such that
for all $x,y\in X\setminus\cals$, $q(x)=q(y)$ if and only if $x,y$ are in the same leaf.
\qed
\end{prop}

Since there are finitely many orbits of singular leaves in $X$,
we can choose such a set $\cals$ containing every singular leaf.

Assume that $Y_v\cap Y_{w}$ is non-degenerate. 
There is  an uncountable number of regular leaves  of $X$  meeting both $\ol U_v\cap \ol U_{w}$.
Since $\ol U_v\setminus U_v\subset \Tilde C$ is contained in a union of singular leaves,
there is a regular leaf meeting both $U_v$ and $U_w$, so $U_v=U_w$, and $v=w$.

We denote by $\Xt$ the subset of the $1$-skeleton of $X$ consisting of the union of all closed transverse edges.
All the paths we consider are chosen as a concatenation of leaf segments and of arcs in $\Xt$.
Let $I$ be an arc in $T$. 
Let $a,b\in\Sigma$ be two preimages of the endpoints of $I$ in $\Xt$.
Choose a path $\gamma$ in $X$ joining $a$ to $b$. 
We view such a path both as a subset of $X$ and as a map $[0,1]\ra X$.
Since $T$ is an $\bbR$-tree, $q(\gamma)\supset I$.
Since the cut locus $C$ is a finite union of leaf segments, 
$\gamma\cap\Tilde C$ consists of a finite number of points and of a finite number
of leaf segments.
In particular, $\gamma$ is a concatenation of a finite number of paths which are contained in $\Tilde C$ or in some $\ol U_v$.
Since a path contained in $\Tilde C$ is mapped to a point in $T$, 
we get that $I$ is covered by finitely many elements of $\caly$.

To prove that the subtree $Y_v$ is closed in $T$, it is sufficient 
to check that given a semi-open interval $[a_0,b_0)$ contained in $Y_v$,
$b_0$ is contained in $Y_v$. Consider $ a, b\in \Xt$ some preimages of $a_0,b_0$, and a 
path $\gamma:[0,1]\ra X$ joining $a$ to $b$. 
By restricting $\gamma$ to a smaller interval, we can assume that for all $t<1$,
$q\circ \gamma(t)\neq b_0$.
Since the image of each transverse edge of $X$ is an arc in $T$,
$K=q(\gamma)$ is a finite tree, and $b_0$ is a terminal vertex of $K$.
In particular, there exists $\eps>0$ such that for $q([1-\eps,1))\subset [a_0,b_0) \subset Y_v$.
We can also assume that $\gamma([1-\eps,1])$ is contained in a transverse edge of $X$. 
For all but countably many $t$'s, $\gamma(t)\notin \cals$ and lies in a regular leaf.
Since $q\circ\gamma(t)\in Y_v$, the leaf $l$ through $\gamma(t)$ meets $\ol U_v$. 
Since $\Tilde C$ consists of singular leaves, $l\subset U_v$.
In particular, $\gamma(t)\in U_v$, so $\gamma(1)\in\ol U_v$ and $b_0\in Y_v$.
We have proved that $\caly$ is a transverse covering of $T$.

Assume that $v$ is such that any leaf of $\pi(U_v)$ is compact.
By definition of the cut locus, this means that every leaf of $U_v$ is regular.
In particular, $U_v$ is a union of leaves of $X$ and $q(U_v)$ is open in $T$.
Since every transverse edge of $X$ is embedded into $T$,
the holonomy along any regular leaf is trivial.
It follows that $U_v$ is a foliated product:
is homeomorphic to $U_v\simeq J\times l$ foliated by $\{*\}\times l$,
where $J$ is an open interval in a transverse edge.
Therefore, $q(U_v)$ is an open interval isometric to $J$.
Since $q(U_v)$ is open in $T$, $Y_v$ is an arc containing no branch point of $T$, 
except possibly at its endpoints.


Assume that  leaves of $\pi(U_v)$ are dense.
We shall prove that $Y_v$ is indecomposable.
First, we claim that any point in $Y_v$ has a preimage in $U_v$.
This will follow from the fact that for any $x\in\ol U_v\setminus U_v$, the leaf through $x$
intersects $U_v$. 
Assume on the contrary that there is a leaf $l$ of $X$ which meets $\ol U_v$ but does not intersect $U_v$.
Consider $l_0$ a connected component of $l\cap\ol U_v$. This is a leaf of the foliated $2$-complex $\ol U_v$.
Let $G_v$ be the global stabilizer of $U_v$, and consider the foliated $2$-complex $\ol U_v/G_v$.
The natural map $\ol U_v/G_v \ra \Sigma$
restricts to an isomorphism between $U_v/G_v$ and a connected component  $\Sigma\setminus C$ (namely $\pi(U_v)$).
Since $U_v/G_v$ consists of finitely many open cells, 
$\ol U_v/G_v$ is compact.
Moreover $\ol U_v/G_v\setminus U_v/G_v$ consists of finitely many vertices and finitely many edges contained in leaves.
The image $\lambda_0$ of $l_0$ in $\ol U_v/G_v$ is a leaf which does not intersect $U_v/G_v$.
Therefore, this leaf is compact
and so is any leaf of $\ol U_v/G_v$ close to $\lambda_0$.
This contradicts the fact that every leaf of $\pi(U_v)$ is dense.

Let's prove indecomposability. Let $I_0,J_0$ be two arcs in $Y_v$.
By taking $I_0$ smaller, one can assume that $I_0=q(I)$ for some arc $I\subset U_v$ contained in a transverse edge.
Consider $a,b$ some preimages of the endpoints of $J_0$ in $U_v\cap \Xt$.
Consider a path $\gamma\subset U_v$ joining $a$ to $b$, and write $\gamma$ as a concatenation 
$l_0.\tau_1.l_1.\tau_2\dots \tau_p.l_p$ where $\tau_i$ is contained in a transverse edge and $l_i$ is a leaf segment.
Since every leaf is dense in $\pi(U_v)$, for all $x\in \tau_1\cup\dots \cup \tau_p$, there exists $g_x\in G$ and a leaf segment $l_x$ in $U_v$
joining $x$ to $g_x.\rond I$. Since any leaf segment in $U_v$ is regular, this is still valid with the same $g_x$ on a neighbourhood of
$x$. By compactness, one can choose each $g_x$ in a finite set $F=\{g_1,\dots,g_n\}$. 
Clearly, $J_0$ is contained in $g_1.q(I)\cup\dots\cup g_n.q(I)$.
Denote by $x$ and $y$ the endpoints of $l_i$.
There is a (regular) leaf segment in $U_v$ joining $g_x.\rond I$ to $g_y.\rond I$,
so $q(g_x.I)\cap q(g_y.I)$ is non-degenerate.
Indecomposability follows,
which completes the proof of Proposition \ref{prop_decompo}.
\end{proof}

\begin{rem}\label{rem_decompo_geom}
  Each indecomposable vertex action $Y_v$ is geometric, and dual to the $2$-complex $\ol U_v$.
Moreover, if $\pi_1(X)$ is generated by free homotopy classes of curves contained in leaves,
then so is $\pi_1(\ol U_v)$ (this will be useful in Appendix \ref{sec_Sela}).
Indeed, consider a curve $\gamma\subset \ol U_v$, and a map $f$ from a planar surface $S$ to $X$
such that one boundary component is mapped to $\gamma$, and every other boundary component
is mapped to a leaf.
Let $S_0\subset S$ be the connected component of $f\m(\ol U_v)$ containing $\ol U_v$.
Since every connected component of $\partial \ol U_v$ is contained in a leaf,
each component of $\partial S_0$ is mapped into a leaf, which proves that $\pi_1(\ol U_v)$ 
is generated by free homotopy classes of curves contained in leaves.
\end{rem}

\subsubsection*{Strong approximation by geometric actions}

We now quote a result from \cite{LP}.
\begin{thm}[{\cite[Theorem 3.7]{LP}}]\label{thm_strongCV}
  Consider a minimal action of a finitely generated group $G$ on an $\bbR$-tree $T$.

Then $G\actson T$ is a strong limit of a direct system of geometric actions $\{(\Phi_k,F_k):G_k\actson T_k\ra G\actson T\}$ 
such that
\begin{itemize*}
\item $\Phi_k$ is one-to-one in restriction to each arc stabilizer of $T_k$,
\item $T_k$ is dual to a $2$-complex $X$ whose fundamental group is generated
by free homotopy classes of curves contained in leaves.
\end{itemize*}\qed
\end{thm}

\begin{rem*}
  The second statement follows from the construction in \cite[Theorem 2.2]{LP}.
\end{rem*}

Since $T_k$ is geometric, it has a decomposition into a graph of actions as described above.
We show how to adapt this result to actions of finitely generated pairs.

\begin{prop}\label{prop_strongCV_rel}
  Consider a minimal action $(G,\calh)\actson T$ of a finitely generated pair on an $\bbR$-tree.
Then there exists a direct system of minimal actions $(G_k,\calh_k)\actson T_k$ converging strongly to $G\actson T$ 
and such that
\begin{itemize*}
\item $\phi_k$ and $\Phi_k$ are one-to-one in restriction to arc stabilizers of $T_k$
\item  $\phi_k$ (resp. $\Phi_k$) restricts to an isomorphism between $\calh_k$
and $\calh_{k+1}$ (resp. $\calh$).
\item $T_k$ splits as a graph of actions where each non-degenerate vertex action is indecomposable, 
or is an arc containing no branch point of $T_k$ except at its endpoints.
\end{itemize*}
\end{prop}

\begin{proof}
We follow Theorem 2.2 and 3.5 from \cite{LP}.
Let $S\subset G$ be a finite set such that $S_\infty=S\cup H_1\dots\cup H_p$ generates $G$.
Let $S_k$ be an exhaustion of $S_\infty$ by finite subsets.
Let $K_k$ be an exhaustion of $T$ by finite trees.
We assume that for each $i\in\{1,\dots,p\}$, $K_k$ contains a point $a_i$ fixed by $H_i$.

 For each $s\in S_k$, consider
$A_s=K_k\cap s\m K_k$ and $B_s=s.A_s=K_k\cap s.K_k$.
We obtain a foliated $2$-complex $\Sigma_k$ as follows.
For each $s\in S_k$, consider a band $A_s\times [0,1]$ foliated by $\{*\}\times [0,1]$,
glue $A_s\times \{0\}$ on $K$ using the identity map on $A_s$,
and glue $A_s\times \{1\}$ using the restriction of $s:A_s\ra K_k$.
For each $s\in H_i\setminus S_k$, add an edge and glue its endpoints on $a_i$;
this edge is contained in a leaf.

The fundamental group of $\Sigma_K$ is isomorphic to the free group $F(S_\infty)$.
Let $\phi:F(S_\infty)\ra  G$ be the natural morphism.
Let $N_k$ be the subgroup of $\ker\phi$ generated by free homotopy classes of curves
contained in leaves (in particular, any relation among the elements of some $H_i$ lies in $N_k$).
This is a normal subgroup of $F(S_\infty)$, and let $G_k=F(S_\infty)/N_k$.
We denote by $\phi_k:G_k\ra G_{k+1}$ and $\Phi_k:G_k\ra G$ the natural morphisms.
Clearly, $G$ is the direct limit of $G_k$.
Let $H_i^{(k)}$ be the image in $G_k$ of the subgroup of $F(S_\infty)$ generated by the elements of $H_i$.
Let $\calh_k=\{H_1^{(k)},\dots,H_p^{(k)}\}$.
The pair $(G_k,\calh_k)$ is finitely generated, and $\phi_k$ and $\Phi_k$ restrict to isomorphisms
between $H_i^{(k)}$, $H_i^{(k+1)}$, and $H_i$.

Let $\Tilde \Sigma_k$ be the Galois covering of $\Sigma_k$ corresponding to $N_k$.
By definition of $N_k$, $\pi_1(\Tilde \Sigma_k)$ is generated by free homotopy classes of curves
contained in leaves.
Let $T_k$ be the leaf space made Hausdorff of $\Tilde \Sigma_k$.
This is an $\bbR$-tree by Theorem \ref{thm_LP_leaves}.
Each $H_i^{(k)}$ fixes a point in $T_k$.
By Proposition \ref{prop_LP_separation}, each arc stabilizer of $T_k$ embeds into $G$.
Finally, the argument of \cite[Th. 2.2]{LP} applies to prove that $G_k\actson T_k$
converges strongly to $G\actson T$.
Minimality follows from Lemma \ref{lem_strong_minimal}.

The $2$-complex $\Sigma_k$ is a finite foliated complex, with a finite set of infinite roses
attached. In particular, $\Sigma_k$ has the same dynamical decomposition as a finite foliated $2$-complex
so one can repeat the argument of Proposition \ref{prop_decompo}.
\end{proof}

\section{Extended Scott's Lemma}\label{sec_scott}

Scott's Lemma claims that if $G$ is a direct limit of groups $G_k$
having compatible decompositions into free products, then $G$ itself has such a decomposition 
(this follows from \cite[Th 1.7]{Scott_coherent}).
Scott's Lemma is usually proved using Scott's complexity (\cite[Th 1.7]{Scott_coherent}).
Delzant has defined a refinement of this complexity for morphisms 
which has many important applications (\cite{Swarup_Delzant}).

The main result of this section is an extension of Scott's Lemma
for more general splittings.
This will be an essential tool to prove piecewise stability of $T$.

In the following statement, an \emph{epimorphism} $(\phi_k,f_k)$ consists of
an onto morphism $\phi_k:G_k\onto G_{k+1}$, and of a continuous map $f_k$ sending an edge to
a (maybe degenerate) edge path.

\begin{SauveCompteurs}{scott}
\begin{thm}[Extended Scott's Lemma]\label{thm_scott}
Let $G_k\actson S_k$ be a sequence of non-trivial actions of finitely generated groups on simplicial trees,
and $(\phi_k,f_k):G_k\actson S_k\ra G_{k+1}\actson S_{k+1}$ be epimorphisms.
Assume that $(\phi_k,f_k)$ does not increase edge stabilizers in the following sense:

$$\forall e\in E(S_k),\forall e'\in E(S_{k+1}),\quad e'\subset f_k(e)\Rightarrow G_{k+1}(e')=\phi_k(G_k(e))\quad (*)$$

Then $\displaystyle\lim_{\ra} G_k$ has a non-trivial splitting over 
the image of an edge stabilizer of some $S_k$.
\end{thm}
\end{SauveCompteurs}

\begin{rem*}
  When edge stabilizers of $S_k$ are trivial, we obtain Scott's Lemma.
\end{rem*}

\subsection{Decomposition into folds.}

\paragraph{Collapses, folds, and group-folds.}

We recall and adapt definitions of \cite{Dun_folding,BF_complexity}.
Consider a finitely generated group $G$ acting on a simplicial tree $S$, without inversion (no element of $G$ flips an edge).
Given an edge $e$ of $S$, collapsing all the edges in the orbit of $e$ 
defines a new tree $S'$ with an action of $G$. We say that $S'$ is a \emph{collapse} of $S$ and we call the
natural map $S\ra S'$ a \emph{collapse}.

Consider two distinct oriented edges $e_1=uv_1$, $e_2=uv_2$ of $S$ having the same origin $u$;
assume that $uv_1$ and $v_2u$ are not in the same orbit (as oriented edges).
Identifying $g.e_1$ with $g.e_2$ for every $g\in G$ defines a new tree $S'$ on which $G$ acts
(without inversion).
We say that $S'$ is obtained by \emph{folding $e_1$ on $e_2$} and
we call the natural map $S\ra S'$  a \emph{fold}.

Consider a vertex $v\in S$, and $N_v$ a normal subgroup of $G(v)$.
Let $N$ be the normal subgroup of $G$ generated by $N_v$ and let $G'=G/N$.
The graph $S'=S/N$ is a tree (see fact \ref{fact_elliptic} below) on which $G'$ acts.
We say that $S'$ is obtained from $S$ by a \emph{group-fold} and we
call the natural map $S\ra S'$ a \emph{group-fold}.
This map is $\phi$-equivariant where $\phi:G\ra G'$ is the quotient map. 

\paragraph{Decomposition into folds.}

The following result is a slight variation on a result by Dunwoody (\cite[Theorem 2.1]{Dun_folding}) . 
For previous results of this nature, see \cite{BF_complexity,Stallings_topology}.

\begin{prop}\label{prop_fold}
Let $(\Phi,F):G\actson T\ra G'\actson T'$ be an epimorphism between finitely supported simplicial actions.
Assume that  
$$\forall e\in E(T),\ \forall e'\in E(T'),\quad e'\subset F(e)\Rightarrow G'(e')=\Phi(G(e)) \quad(*)$$

Then we may subdivide $T$ and $T'$ 
so that there exists a finite sequence of simplicial actions $G\actson T=G_0\actson T_0,\dots,G_n\actson T_n=G'\actson T'$, 
 some epimorphisms $(\phi_i,f_i):G_i\actson T_i\ra G_{i+1}\actson T_{i+1}$ such that $\Phi=\phi_{n-1}\circ\dots\circ \phi_0$,
each $(\phi_i,f_i)$ satisfies $(*)$, and is either 
\begin{enumerate*}
\item a collapse
\item a group-fold
\item or a fold between two edges 
$uv,uv'$ of $T_i$ such that $G_i(u)$ injects into $G'$ under $\Phi_{i}$.
\end{enumerate*}
\end{prop}

\begin{rem*}
  We don't claim that $F=f_{n-1}\circ\dots\circ f_1$.
\end{rem*}

\begin{proof}
We may change $F$ so that it is linear in restriction to each edge of $T$.
Then we can subdivide $T$  and $T'$ so that $F$ maps each vertex to a vertex and each edge to an edge or a vertex.
The new map satisfies $(*)$.
Set $T_0=T$, $G_0=G$.

We describe an iterative procedure.
We assume that the construction has begun, and that we are given
$(\Phi_i,F_i):G_i\actson T_i\ra G'\actson T'$ satisfying $(*)$.

Step 1. Assume that $F_i$ maps an edge $e$ of $T_i$ to a point.
Let $T_{i+1}$ be the tree obtained by collapsing $e$ and let $G_{i+1}=G_i$, $\phi_i=\id$ and $\Phi_{i+1}=\Phi_i$.
Define $f_i:T_i\ra T_{i+1}$ as the collapse map, and $F_{i+1}:T_{i+1}\ra T'$ as the map induced by $F_i$.
Then return to step 1.

Clearly, $(\id,f_i)$ and $(\Phi_{i+1},F_{i+1})$ satisfy $(*)$.
Since $T_{i+1}$ has fewer orbits of edges than $T_i$, we can repeat step 1 until
$F_i$ does not collapse any edge of $T_i$.

Step 2.
Assume that there is a vertex $v\in T_i$ such that the kernel $N_v$ of $\Phi_i{}_{|G_i(v)}:G_i(v)\ra G'$ is non-trivial.
Let $N$ the normal closure of $N_v$ in $G_i$, and consider $G_{i+1}=G_i/N\actson T_{i+1}=T_i/N$
the action obtained by the corresponding group-fold.
Define $f_i:T_i\ra T_{i+1}$ and $\phi_i:G_i\ra G_{i+1}$ as the quotient maps.
Define $F_{i+1}:T_{i+1}\ra T'$ and $\Phi_{i+1}:G_{i+1}\ra G'$ as the maps induced by $F_i$ and $\Phi_i$.
Then return to step 2.

The group fold $(\phi_i,f_i)$ automatically satisfies (*).
Moreover,  $(\Phi_{i+1},F_{i+1})$ inherits property (*) from $(\Phi_i,F_i)$.
This step will repeated only a finite number of times since step 2 decreases the number of orbits
of vertices with $N_v\neq\{1\}$.

Step 3.
Assume that there exists two edges $e_1=uv_1$, $e_2=uv_2$ of $T_i$ such that $F_i(e_1)=F_i(e_2)$.
Define $T_{i+1}$ as the tree obtained by folding $e_1$ on $e_2$,  $f_i:T_i\ra T_{i+1}$ as the folding map,
and $F_{i+1}:T_{i+1}\ra T'$ as the induced map. Define $G_{i+1}=G_i$, $\phi_i=\id$ and $\Phi_{i+1}=\Phi_i$.
Then return to step 2. 

Denote by $\eps$ the common image of $e_1$ and $e_2$ in $T_{i+1}$.
Since at the beginning of step 3, no group fold can be done,
$\Phi_i{}_{|G_i(u)}$ is one-to-one. 
We claim that $e_1$ and $e_2$ cannot be in the same orbit.
Indeed, if $e_1=g.e_2$ for some $g\in G$, then $g\in G_i(u)$ (because $T'$ has no inversion)
so $\Phi_i(G_i(e_1))\subsetneq \Phi_i(\langle G_i(e_1),g\rangle)= \Phi_{i+1}(G_{i+1}(\eps))\subset G'(F_{i}(e_1))$,
contradicting $(*)$.
It follows that step 3 will be repeated only finitely many times
because it decreases the number of orbits of edges, and step 2 does not change it.
Since $F_i$ satisfies (*) and $\Phi_i{}_{|G_i(u)}$ is one-to-one, 
we get $\Phi_i(G_i(e_1))=\Phi_i(G_i(e_2))$ so $G_i(\eps)=G_i(e_1)$ and $f_i$ satisfies $(*)$.
The fact that $F_{i+1}$ inherits $(*)$ is clear.

When step 3 cannot be repeated any more, $F_i$
is an isometry. If follows that any $g\in\ker \Phi_i$ fixes $T_i$ pointwise.
Since $\Phi_i{}_{|G_i(u)}$ is one-to-one for every vertex, $\Phi_i$ is an isomorphism.
\end{proof}

\subsection{Proof of Extended Scott's Lemma}

\newcommand{\Ell}{\mathrm{Ell}}
\begin{proof}
One can assume that each $S_k$ is minimal.
By Proposition \ref{prop_fold}, we may assume that each map $f_k$ is a collapse, a group-fold or a fold.
Let $\Ell(S_k)$ be the subset of $G_k$ consisting of elements fixing a point in $S_k$.
The first Betti number of the graph $S_k/G_k$ 
coincides with the rank of the free group $G_k/\langle\Ell(S_k)\rangle$.
In particular, this Betti number is non-increasing and
we can assume that it is constant.

We work at the level of quotient graph of groups $\ol S_k=S_k/G_k$ and we denote by $x\mapsto \ol x$ the quotient map. 
Consider an oriented edge $\ol e$ of $\ol S_k$ with terminal vertex $\ol v=t(\ol e)$.
Say that $\ol e$ carries the symbol $=$
if the edge morphism $i_{\ol e}:G_{\ol e}\ra G_{\ol v}$ is onto.
Otherwise, we say that $\ol e$ carries $\neq$.
Define $W_k$ as the set of vertices $\ol v\in\ol S_k$ such that there is an oriented edge $\ol e$ with $t(\ol e)=\ol v$
carrying $\neq$.
At the level of the tree, $W_k$ is the set of orbits of vertices $v\in S_k$ for which there exists an edge $e$ incident on $v$
with $G_k(e)\subsetneq G_k(v)$. 
Note that $\#W_k$ is invariant under subdivision.

We claim that $\#W_k$ is non-increasing, and that
folds and collapses which don't decrease $\#W_k$ don't change $\Ell(S_k)$.

Assume that $f_k:S_k\ra S_{k+1}$ is induced by the collapse of an edge $e=uv$.
\newcommand{\olf}{\Bar f}
Let $\olf_k:\ol S_k\ra \ol S_{k+1}$ be the induced map.
The endpoints of $\ol e$ are distinct since otherwise, the first Betti number would decrease.
If both $\ol u$ and $\ol v$ are outside $W_k$, then the stabilizer of $u$ coincides with that of $f_k(u)$;
in particular, $\olf_k(\ol u)\notin W_k$ and $\#W_k$ is non-increasing.
If $\Ell(S_k)$ increases under the collapse, then both orientations of $\ol e$ carry $\neq$,
in which case both $\ol u$ and $\ol v$ belong to $W_k$, so $\#W_{k+1}< \# W_k$.

Now assume that $f_k:S_k\ra S_{k+1}$ is the fold of $e_1=uv_1$ with $e_2=uv_2$ 
(in particular $G_k=G_{k+1}$).
Denote by $v'$ (resp $e'$) the common image of $v_1,v_2$ (resp. $e_1,e_2$) in $S_{k+1}$,
and $u'$ the image of $u$.
Since $f_k$ satisfies $(*)$, $e_1$ and $e_2$ are in distinct orbits and $G_k(e')=G_k(e_1)=G_k(e_2)$.
Since the first Betti number of $\ol S_k$ is constant, $\ol v_1\neq \ol v_2$.
It may happen that $\ol u=\ol v_i$ for some $i\in\{1,2\}$.
By subdividing $e_1$ and $e_2$ and replacing the original fold by two consecutive folds, 
we can ignore this case. Thus we assume that $\ol u,\ol v_1,\ol v_2$ are distinct.
One has $G_k(v')=\langle G_k(v_1),G_k(v_2)\rangle$.
It follows that if both $\ol v_1$ and $\ol v_2$ lie outside $W_k$, then $\ol v\notin W_{k+1}$.
Since $G_k(u)=G_k(u')$, $\ol u\in W_k$ if and only if $\ol u'\in W_{k+1}$.
It follows that $\#W_k$ is non-increasing.
If $\Ell(S_k)$ increases, then $G_k(v)=\langle G_k(v_1),G_k(v_2)\rangle$ is distinct from both $G_k(v_1)$
and $G_k(v_2)$. Therefore, $\ol e_1$ can't carry $=$ at $\ol v_1$, since this would imply $G_k(v_1)=G_k(e_1)=G_k(e_2)\subset G_k(v_2)$.
Similarly, $\ol e_2$ carries $\neq$ at $\ol v_2$.
Therefore, $\ol v_1,\ol v_2\in W_k$, and $\#W_{k+1}< \#W_k$.
This proves the claim.
\\

Without loss of generality, we assume that $\#W_k$ is constant. 
For each $k$, either $f_k$ is a group-fold, or $\Ell(S_k)=\Ell(S_{k+1})$.
Let $N_k=\ker \phi_{k-1}\circ\dots\circ\phi_1$.
We prove by induction on $k$ that $S_0/N_k$ is a tree and $\Ell(S_0/N_k)=\Ell(S_k)\subset G_k$.
We will use the following standard fact:
\begin{fact}\label{fact_elliptic}
  Let $T$ be a simplicial tree and $N$ a group of isometries. Then $T/N$
is a tree if and only if $N$ is generated by elliptic elements.
\end{fact}

\begin{proof}
If the graph of groups $T/N$ is a tree, then its fundamental group is
generated by its vertex groups.
If $T/N$ is not a tree, killing its vertex groups gives a non-trivial free group.
\end{proof}

Assume that $S_0/N_k$ is a tree. If $f_k$ is not a group-fold,
then $G_k=G_{k+1}$, $N_k=N_{k+1}$ and $\Ell(S_k)=\Ell(S_{k+1})$ so we are done.
If $f_k$ is a group fold, then $G_{k+1}$ is the quotient of $G_k$ by a normal subgroup $K$ generated by subset of 
$\Ell(S_k)$ and $S_{k+1}=S_k/K$. 
It follows that $S_0/N_{k+1}=(S_0/N_k)/K$ is a tree.
For $g\in G_{k}$, denote by $\ol g$ its image in $G_{k+1}=G_k/K$.
We get: $\ol g\in \Ell(S_k/K)\Leftrightarrow \exists k\in K, gk\in \Ell(S_K)$
$\Leftrightarrow \exists k\in K, gk\in \Ell(S_0/N_k)$
$\Leftrightarrow \ol g\in \Ell(S_0/N_{k+1})$. The induction follows.

Since $S_0/N_k$ is a tree for all $k$, $N_k$ is generated by elliptic elements,
and so is $N=\cup_k N_k$. This way, we get an action of $G=\displaystyle\lim_{\ra} G_k=G_0/N$ on the tree $S_0/N$.
Assume that this action has a global fix point and argue towards a contradiction. 
Let $\{g_1,\dots,g_p\}$ be a generating set of $G_0$.
There exists $x\in S_0$ and $n_1,\dots,n_p\in N$
such that $g_in_i.x=x$. Choose $k$ large enough so that $n_1,\dots,n_p\in N_k$.
Then $G_k$ has a global fix point in $S_0/N_k$. Since  $\Ell(S_0/N_k)=\Ell(S_k)$,
$G_k$ fixes a point in $S_k$, a contradiction.

Finally,  an edge stabilizer in $G\actson S_0/N$ is the image of an edge stabilizer in $G_0\actson S_0$,
which concludes the proof of Extended Scott's Lemma.
\end{proof}

\subsection{Relative version of Extended Scott's Lemma}

\begin{thm}[Extended Scott's Lemma, relative version.]  \label{thm_scott_rel}
Consider $(G_k,\{H_1^k,\dots H_p^k\})\actson S_k$ a sequence of non-trivial actions of finitely generated pairs on simplicial trees.
Let $(\phi_k,f_k):G_k\actson S_k\ra G_{k+1}\actson S_{k+1}$ be epimorphisms
mapping $H_i^k$ onto $H_i^{k+1}$.
Consider $G=\displaystyle\lim_{\ra} G_k$ the inductive limit and $\Phi_k:G_k\ra G$ the natural map.

Assume that
\begin{enumerate*}
\item\label{hyp1} $\forall e\in E(S_k),\forall e'\in E(S_{k+1}),\quad e'\subset f_k(e)\Rightarrow G_{k+1}(e')=\phi_k(G_k(e))\quad (*)$
\item\label{hyp2} $\forall e\in E(S_k),\forall i$, $\Phi_k(H_i^k)\not\subset\Phi_k(e)$ 
\end{enumerate*}

Then the pair $(G,\calh)$  has a non-trivial splitting over 
the image of an edge stabilizer of some $S_k$.
Moreover, any subgroup $H\subset G$ fixing a point in some $S_k$ fixes a point in the obtained splitting of $G$.
\end{thm}

\begin{rem}
  The additional assumption is necessary 
as shows the following example.
Let $A$ be an unsplittable finitely generated group containing a finitely generated free group $F$.
Let $\{a_1,a_2,\dots\}\subset F$ be an infinite basis of a free subgroup of $F$,
and $F'=F(a'_1,a'_2,\dots)$ another free group with infinite basis.
Consider $$G_k=\left\langle F*F'\ |\  a_1=a'_1, \dots, a_k=a'_k\right\rangle*_F A.$$
The group $G_k$ is finitely generated relative to $F'$.
This sequence of splittings (over $F$) satisfies condition $(*)$ but not the additional assumption.
The inductive limit of $G_k$ is the unsplittable $A$.
\end{rem}

\begin{proof}
Since an action of a finitely generated pair is finitely supported,
we can use the decomposition into folds of Proposition \ref{prop_fold}.
The proof  of Extended Scott's Lemma does not use finite generation until 
the proof of the non-triviality of the obtained splitting.
Recall the notations of the end of the proof of Theorem \ref{thm_scott}:
$\Phi_k:G_k\ra G$ is the natural morphism,
$N=\ker \Phi_0:G_0\ra G$, $N_k=\ker\phi_{0k}:G_0\ra G_k$ so that $N=\bigcup_k N_k$.

At this point, we know that we can forget finitely many terms in our sequence of actions
so that $S_0/N$ is a tree, and that for all $k$, the action $G_k\actson S_0/N_k$ is non-trivial.
We assume that $G\actson G_0/N$ has a global fix point and argue towards a contradiction.
Let $S^0\subset G_0$ be a finite set such that $S^0\cup H^0_1\cup\dots\cup H_p^0$ generates $G$.
Consider $\{g_1,\dots,g_q\}$ such that  $\{g_1,\dots,g_q\}\cup H^0_1\cup\dots\cup H_p^0$ generates $G_0$.
There exists $a\in S_0$ and $n_1,\dots,n_q\in N$
such that $g_jn_j.a=a$ for all $j\in\{1,\dots,q\}$. 
Let $b_i\in S_0$ be a fix point of $H_i^0$.
By the second hypothesis, $H_i^0$ fixes no edge in $S_0/N$ so it fixes a unique point in $S_0/N$.
In particular, the images of $b_i$ and $a$  in $S_0/N$ coincide. 
Therefore, $b_i=n'_i.a$ for some $n'_i\in N$.
Choose $k$ large enough so that $N_k$ contains all those elements $n_j,n'_i$.
Then $G_k$ fixes the image of $a$ in $S_0/N_k$.
This contradicts the non-triviality of $G_k\actson S_0/N_k$.
\end{proof}

\section{Getting piecewise stability}\label{sec_acc2pw}

\begin{thm}\label{thm_acc2pw}
Consider a minimal action of a finitely generated pair $(G,\calh)$ on an $\bbR$-tree $T$.
Assume that
\begin{enumerate*}
\item $T$ satisfies the ascending chain condition;
\item there exists a finite family of arcs $I_1,\dots,I_p$ such that $I_1\cup\dots\cup I_p$ spans $T$ 
and such that for any unstable arc $J$ contained in some $I_i$,
\begin{enumerate*}
\item $G(J)$ is finitely generated;
\item  $G(J)$ is not a proper subgroup of any conjugate of itself
\ie $\forall g\in G$, $G(J)^g\subset G(J)\Rightarrow G(J)^g= G(J)$.
\end{enumerate*}
\end{enumerate*}

Then either $(G,\calh)$ splits over the stabilizer of an unstable arc contained in some $I_i$,
or $T$ is piecewise-stable.
\end{thm}

\begin{proof}
By enlarging some peripheral subgroups,
we may assume that each $H\in\calh$ is either finitely generated, 
or is not contained in a finitely generated elliptic subgroup of $G$.

By Proposition \ref{prop_strongCV_rel}, consider $(G_k,\calh_k)\actson T_k$ a sequence of actions on $\bbR$-trees converging strongly
to $T$, and such that $T_k$ splits as a graph of actions where each vertex action is either 
an indecomposable component, or an edge.
Denote by $\phi_k:G_k\ra G_{k+1}$, $\Phi_k:G_k\ra G$,
$f_k:T_k\ra T_{k+1}$, and $F_k:T_k\ra T$ the morphisms of the corresponding direct system.

Without loss of generality, we can assume that there exists
$\Tilde I_1,\dots, \Tilde I_p\subset T_0$ which map isometrically to  $I_1,\dots, I_p$. 
By subdividing the arcs $\Tilde I_i$ and the edges occurring in the decomposition of $T_k$, 
we may assume that each $\Tilde I_i$ is either an edge or
an indecomposable arc of $T_0$.
All the arcs of $T_k$ we are going to consider embed isometrically into $T$.
Say that an arc $I\subset T_k$ is \emph{pre-stable} if its image in $T$ (under $F_k$) is a stable arc.
Otherwise, we say that $I$ is \emph{pre-unstable}.
If $I$ is contained in an indecomposable component, then it is pre-stable by Lemma \ref{lem_indec_stab} 
and Assertion 1 of Lemma \ref{lem_indecompo}.



To do some bookkeeping among pre-unstable edges, we shall construct inductively a combinatorial graph  
$\calt=\calt(\Tilde I_1)\dunion \dots  \dunion\calt(\Tilde I_p)$
as a disjoint union of rooted trees.
We use standard terminology for rooted trees: 
the father of $v$ is the neighbour of $v$ closer to the root than $v$,
a child of $v$ is a neighbour of $v$ which is not the father of $v$, 
an ancestor of $v$ if a vertex on the segment joining $v$ to the root.
The level of a vertex is its distance to the root.
We denote by $\calt_k(\Tilde I_i)$ (resp.\ $\calt_{\geq k}(\Tilde I_i)$)
the set of vertices of $\calt(\Tilde I_i)$ of level $k$ (resp. at least $k$),
and $\calt_k=\calt_k(\Tilde I_1)\cup\dots\cup\calt_k(\Tilde I_p)$.
(resp.\ $\calt_{\geq k}=\calt_{\geq k}(\Tilde I_1)\cup\dots\cup\calt_{\geq k}(\Tilde I_p)$).

Each vertex $v$ of level $k$ of $\calt(\Tilde I_i)$ will be labeled by a pre-unstable edge $J_v$
contained in $f_{0k}(\Tilde I_i)\subset T_k$.
We label the root of $\calt(\Tilde I_i)$ by $\Tilde I_i$.
Assume that $\calt$ has been constructed up to the level $k$.
Subdivide the edge structure of $T_{k+1}$ so that, for each  $v\in\calt_k$,
$f_k(J_v)$ is a finite union of edges and of indecomposable arcs.
We may also assume that the length of each edge is at most $1/2^k$.
The indecomposable pieces of $f_k(J_v)$ are pre-stable, and we discard them.
We also discard pre-stable edges of $f_k(J_v)$.
For each pre-unstable edge $J'$ contained in $f_k(J_v)$, we add new child of $v$ labeled by $J'$.

If $\calt$ is finite,
then each $\Tilde I_i$ is contained in a finite union of stable arcs.
Since $\Tilde I_1\cup\dots \Tilde I_p$ spans $T$, $T$ is piecewise stable. So we assume that $\calt$ is infinite.

To each vertex $v\in\calt_k$, we attach two subgroups of $G$:
$A_v=\Phi_k(G_k(J_v))$, and $B_v=G(F_k(J_v))$. Clearly, $A_v\subset B_v$.
Given $u,v\in\calt$, write $B_u< B_{v}$ if $B_{u}$ is properly contained in some conjugate of $B_{v}$.
Say that $v\in\calt_k$ is \emph{minimal} if there is no $u\in\calt_k$ with $B_u < B_v$.
Let $\calm\subset\calt$ be the set of minimal vertices, and $\calm_k=\calm\cap\calt_k$.
Since no $B_v$ is  a proper conjugate of itself, for each $v$, either $v\in\calm_k$ or there exists $u\in\calm_k$ with $B_u<B_v$.

\begin{lem}\label{lem_immortel}
There exists $k_0$ such that for all $k\geq k_0$, the following hold:
\begin{enumerate*}
\item \label{it_notmin} if $v\in \calt_k\setminus\calm_k$, then no child of $v$ is minimal;
\item \label{it_min} if $v\in \calm_k$ then for any child $v'\in\calm_{k+1}$ of $v$, $A_v=A_{v'}$;
\end{enumerate*}
\end{lem}

\begin{proof}
Say that $v$ is a \emph{clone} of $v'$ if $B_v=B_{v'}$ and $A_v=A_{v'}$.
A \emph{genealogical line} is a  sequence $v_0,v_1,\dots$ of vertices of $\calt$
where $v_{i+1}$ is a child  of $v_i$.
Say that $v\in \calt$ is \emph{immortal} if there is a genealogical line starting at $v$
and consisting of clones of $v$.

We claim that any genealogical line $v_0,v_1,\dots$ eventually consists of clones. 
By the ascending chain condition in $T$, there are at most finitely many indices $i$ such that
$B_{v_i}\subsetneq B_{v_{i+1}}$. So we can assume that $B_{v_i}=B_{v_0}$ for all $i$.
Since $J_{v_0}$ is pre-unstable, $B_{v_0}$ is finitely generated.
The strong convergence implies that for $k$ large enough,
every generator of $B_{v_0}$ has a preimage in
$G_{k}(f_{k_0k}(J_{v_0}))$ (where $k_0$ is the level of $v_0$). 
Since for each $v_i$ of level at least $k$,
$J_{v_i}\subset f_{k_0k}(J_{v_0})$, we get 
$B_{v_0}\subset A_{v_i}\subset B_{v_i}=B_{v_0}$.
Therefore, all the vertices $v_i$ of sufficiently large level are clones of each other.
This proves the claim. In particular, this also proves that if $v$ is immortal, then $A_v=B_v$.

Say that $v$ is \emph{post-immortal} if it is immortal or if it has an immortal ancestor.
We claim that all but finitely many vertices of $\calt$ are post-immortal.
Indeed, for each $i\in\{1,\dots,p\}$, the set of vertices of $\calt(\Tilde I_i)$ which are not post immortal is a 
rooted subtree of $\calt(\Tilde I_i)$. If it is infinite, it contains a genealogical line.
Since this line contains many immortal vertices, this is a contradiction.

Let $k_0$ be such that $\calt_{\geq k_0}$ consists of post-immortal vertices.
We shall use several times the following fact:
any $v\in\calm_{\geq k_0}$ is the clone of an immortal vertex.
Indeed, let $u'$ be an immortal ancestor of $v$, and let $u$ be a clone of
$u'$ of same level as $v$;
then $A_{u'}=B_{u'}\subset A_v\subset B_v$, and since $v$ is minimal, inclusions are equalities,
and $u'$ is a clone of $v$.
It follows that for all $v\in\calm_{\geq k_0}$, $A_v=B_v$.

If $v\in \calt_k\setminus\calm_k$, there exists $u\in \calm_k$ with $B_u<B_v$.
Since $u$ is the clone of an immortal vertex, there exists $u'$ of level $k+1$ with $B_u=B_{u'}$ 
and no child of $v$ can be minimal. This proves assertion \ref{it_notmin}.

If $v\in \calm_k$, and $v'\in \calm_{k+1}$ is a child of $v$,
then $v$ has an immortal clone, so has an immortal clone $u'$ of level $k+1$.
We have $B_{u'}=A_{u'}=A_v\subset A_{v'}\subset B_{v'}$ and
since $v'$ is minimal, $v'$ is a clone of $v$. Assertion \ref{it_min} follows.
\end{proof}

We always consider $k\geq k_0$.
Let $\cale_k$ be the set of edges of our decomposition of $T_k$ which are in the orbit of
some $J_v$ for some $v\in\calm_k$. 
In our decomposition of $T_k$ as a graph of actions, collapse all the vertex trees which are not in $\cale_k$
 (see Definition \ref{dfn_collapse}).
The resulting tree $S_k$ is a graph of actions where each vertex action is an edge, so it is a simplicial tree.
The action $G_k\actson S_k$ is minimal by Lemma \ref{lem_collapse} and is therefore non-trivial.
By assertion \ref{it_notmin} of Lemma \ref{lem_immortel}, 
 $f_k$ induces a natural map $S_k\ra S_{k+1}$.
Since any $H\in \calh_k$ fixes a point in $T_k$, it also does in $S_k$, so the action on $S_k$ is an action of the pair $(G_k,\calh_k)$.

Let's check that the actions $G_k\actson S_k$ satisfy the hypotheses of Extended Scott's lemma.
Consider an edge $e\in\cale_k$. We may translate $e$ so that
$e$ corresponds to an arc $J_v$, with $v\in \calm_k$. 
Let $e'\in \cale_{k+1}$ with $e'\subset f_k(e)$.
Since $e'\in \cale_{k+1}$, it is pre-unstable, so $e'=J_{v'}$ for some child $v'\in\calt$ of $v$.
Moreover, since $e'\in\cale_{k+1}$, $e'$ is in the orbit of some $J_{w'}$ with $w'$ minimal,
so $v'\in \calm_{k+1}$.
By assertion \ref{it_min} of  Lemma \ref{lem_immortel}, $\Phi_{k+1}(G_{k+1}(e'))=A_{v'}=A_v=\Phi_k(G_k(e))$.
Since $\Phi_k$ and $\Phi_{k+1}$ are  one-to-one in restriction to arc stabilizers,
we get $\phi_k(G_k(e))=G_{k+1}(e')$.

If $\calh=\es$, this is enough to apply Extended Scott's Lemma,
so $G$ splits over some group $A_v$.
Since $A_v=B_v$ is the stabilizer of an unstable arc in $T$, the theorem is proved in the non-relative case.

If $\calh\neq\es$, we need to modify slightly the argument to ensure that the second hypothesis
of the relative version of Extended Scott's Lemma holds.
Let $\calh'\subset \calh$ be the subset consisting of subgroups 
which are not finitely generated and $\calh'_k\subset \calh_k$
the subset corresponding to $\calh'$.
The pairs $(G,\calh')$ and $(G_k,\calh'_k)$ are finitely generated.
Recall that no $H\in \calh'$ is contained in a finitely generated elliptic subgroup.
In particular, for all $H\in\calh'_k$, $\Phi_k(H)$ fixes no unstable arc of $T$.
For each $e\in E(S_k)$, $\Phi_k(G_k(e))$ fixes an unstable arc in $T$ so
Theorem \ref{thm_scott_rel} applies.
We thus get a non-trivial splitting  $(G,\calh')\actson S$. 

In view of the moreover part of Theorem \ref{thm_scott_rel},
since any $H\in\calh_0$ fixes a point in $S_0$,
any $H\in\calh$ fixes a point in $S$.  
Thus, $S$ defines a non-trivial splitting of $(G,\calh)$.
\end{proof}

\section{Piecewise stable actions}\label{sec_pw2triv}

The goal of this section is the following result:
\begin{thm}\label{thm_pw2triv}
  Let $(G,\calh)$ be a finitely generated pair having a piecewise stable action
on an $\bbR$-tree $T$. 

Then, either $(G,\calh)$ splits over the stabilizer $H$ of an infinite tripod (and the normalizer of $H$
contains a non-abelian free group generated by two hyperbolic elements whose axes don't intersect), 
or $T$ has a decomposition into a graph of actions where each
vertex action is either 
\begin{enumerate*}
\item simplicial: $G_v\actson Y_v$ is a simplicial action on a simplicial tree;
\item of Seifert type: the vertex action $G_v\actson Y_v$ has kernel $N_v$, and the faithful action $G_v/N_v\actson Y_v$
is dual to an arational measured foliation on a closed $2$-orbifold with boundary;
\item axial: $Y_v$ is a line, and the image of $G_v$ in $\Isom(Y_v)$ is a finitely generated group acting with dense orbits on $Y_v$.
\end{enumerate*}
\end{thm}

The proof relies on the following particular case of a result by Sela:
\begin{UtiliseCompteurs}{thm_sela}
\renewcommand{\thesection}{\Alph{section}}
  \begin{thm}[\cite{Sela_acylindrical}]
  Consider a minimal action of a finitely generated group $G$ on an $\bbR$-tree $T$ with trivial arc stabilizers.

Then, either $G$ is freely decomposable, or $T$ has a decomposition into a graph of actions where each
vertex action is either 
\begin{enumerate*}
\item of surface type: the vertex action $G_v\actson Y_v$ 
is dual to an arational measured foliation on a closed $2$-orbifold with boundary;
\item  axial: $Y_v$ is a line, and the image of $G_v$ in $\Isom(Y_v)$ is a finitely generated group acting with dense orbits on $Y_v$.
\end{enumerate*}
\end{thm}
\end{UtiliseCompteurs}

We shall give a proof of this result in Appendix \ref{sec_Sela}.

\subsection{From piecewise stability to trivial arc stabilizers}

\begin{lem}\label{lem_pw2goa}
  Let $G$ be a  group having a piecewise stable action
on an $\bbR$-tree $T$. 

Then $T$ has a decomposition into a graph of actions such that,
denoting by $N_v$ the kernel of the vertex action $G_v\actson Y_v$,
then $G_v/N_v\actson Y_v$ has trivial arc stabilizers.
\end{lem}

\begin{rem*}
The lemma doesn't assume any finite generation of $G$.
\end{rem*}

\begin{proof}
Remember that a subtree $Y\subset T$ is stable if the stabilizer of any arc $J\subset Y$
fixes $Y$.   Given a stable arc $I\subset T$, consider $Y_I$ the maximal stable subtree containing $I$.
This is a well defined subtree because if two stable subtrees contain $I$, their union is still stable,
and an increasing union of stable subtrees is stable.
Moreover, $Y_I$ is closed in $T$ because the closure of a stable subtree is stable.

Let $(Y_v)_{v\in V}$ be the family of all maximal stable subtrees of $T$.
By piecewise stability, any arc of $T$ is contained in a finite union of them.
By maximality, if $Y_u\cap Y_v$ is non-degenerate, then $Y_u=Y_v$.
Therefore, the family $(Y_v)_{v\in V}$ is a transverse covering of $T$.
Denote by $G_v$ the global stabilizer of $Y_v$, and by $N_v$ its pointwise stabilizer.
Since $Y_v$ is a stable subtree, $G_v/N_v\actson Y_v$ has trivial arc stabilizers.
\end{proof}

\subsection{Relative version of Sela's Theorem}\label{sec_relfg}

Consider a decomposition of $T$ as a graph of actions as in Lemma \ref{lem_pw2goa}.
Consider a vertex group $G_v$ of the corresponding graph of groups $\Gamma$.
Even if $G$ is finitely generated, $G_v$ may fail to be finitely generated.
However, in this case, the peripheral structure $(G_v,\calh_v)$ of $G_v$ in $\Gamma$
is finitely generated  (Lemma \ref{lem_supp_fini} and Lemma \ref{lem_relfg}),
and the action on $T_v$ is an action of the pair $(G_v,\calh_v)$.
If we started with a finitely generated pair $(G,\calh)$, $G_v$ is finitely generated with respect to 
a set $\calh_v$ consisting of the peripheral structure of $G_v$ in $\Gamma$ together with some conjugates of elements of $\calh$ 
(Lemma \ref{lem_relfg_rel}).
We will need a version of Sela's result applying in this context.

\begin{prop}[Relative version of Sela's result]\label{prop_sela_rel}
  Consider a minimal action of a finitely generated pair $(G,\calh)$ on an $\bbR$-tree $T$ with trivial arc stabilizers.

Then, either $(G,\calh)$ is freely decomposable, or $T$ has a decomposition into a graph of actions where each
non-degenerate vertex action is either 
\begin{enumerate*}
\item of surface type: the vertex action $G_v\actson Y_v$ 
is dual to an arational measured foliation on a closed $2$-orbifold with boundary;
\item axial: $Y_v$ is a line, and the image of $G_v$ in $\Isom(Y_v)$ is a finitely generated group acting with dense orbits on $Y_v$.
\end{enumerate*}
\end{prop}

\begin{proof}
We shall embed $G$ into a finitely generated group $\Hat G$ and apply Sela's non-relative result.
Let $\calh=\{H_1,\dots,H_p\}$. We may assume that each $H_i$ is non-trivial.
For $i=1,\dots,p$, choose some finitely generated group $\Hat H_i$ containing $H_i$.
We may assume that $\Hat H_i$ is freely indecomposable by changing $\Hat H_i$ to $\Hat H_i\times \bbZ/2\bbZ$.
Consider the graph of groups $\Gamma$ below.
\begin{center}\includegraphics{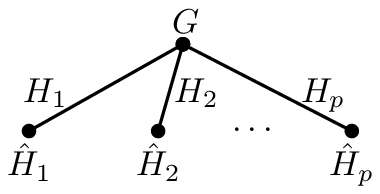}\end{center}
%
We denote by $u$ (resp. $u_i$) the vertex labeled by $G$ (resp. $\Hat H_i$).
 Clearly, the group $\Hat G=\pi_1(\Gamma)$ is finitely generated.

 We shall add a structure of a graph of actions on $\bbR$-trees on $\Gamma$. 
 Let $Y_u$ be a copy of $T$ endowed with its natural action of $G$,
 and let $Y_{u_i}$ be a point endowed with the trivial action of $\Hat H_i$.
 For each edge $uu_i$ of $\Gamma$, we define its attaching point in $Y_u$ 
 as the (unique) point of $T$ fixed by $H_i$.
 Let $\Hat G\actson \Hat T$ be the $\bbR$-tree dual to this graph of actions.

 Since $T$ is a union of axes of elements of $G$,
 and since $\Hat T$ is covered by translates of $T$, 
 $\Hat T$ is a union of axes of elements of $\Hat G$ so $\Hat T$ is minimal.

This action has trivial arc stabilizers because $G\actson T$ does.
Therefore, we can apply Sela's Theorem \ref{thm_Sela} to  $\Hat G\actson \Hat T$.

Assume first that $\Hat G$ is freely decomposable, \ie that $\Hat G$ acts non-trivially on a simplicial
tree $R$ with trivial edge stabilizers. Then $\Hat H_i$ is elliptic in $R$ because it is freely indecomposable.
If $G$ fixes a point $x\in R$, then each $H_i$ fixes $x$ and cannot fix any other point because 
edge stabilizers are trivial. It follows that each $\Hat H_i$ fixes $x$ and that $\Hat G$ fixes $x$, a contradiction.
Thus, the action of $G$ on $R$ defines a non-trivial free decomposition of $G$ relative to $H_i$.

Assume now that $\Hat T$ has a decomposition into
a graph of actions where each vertex action is axial or of surface type.
In particular, each vertex action is indecomposable.
Let $(Z_i)$ be the transverse covering of $\Hat T$ induced by this decomposition.
Since $Z_i$ is indecomposable, if $Z_i\cap T$ is non-degenerate, then $Z_i\subset T$ (Lemma \ref{lem_indec_component}).
The family of subtrees $Z_i$ contained in $T$ is therefore a transverse covering of $T$.
This gives a decomposition of $G\actson T$ as a graph of actions where each vertex action
is axial or of surface type.
\end{proof}

\subsection{Proof  of Theorem \ref{thm_pw2triv}}

\begin{proof}[Proof of Theorem \ref{thm_pw2triv}]
Consider a transverse covering $\caly=(Y_v)_{v\in V}$ of $T$ 
as in Lemma \ref{lem_pw2goa}.
Denote by $G_v$ the global stabilizer of $Y_v$, 
and by $N_v\normal G_v$ the kernel of the vertex action $G_v\actson Y_v$.
We denote by $\ol G_v=G_v/N_v$ the group acting with trivial arc stabilizers on $Y_v$.

We shall prove that either $(G,\calh)$ splits over a tripod stabilizer with the desired properties,
or that each vertex action $G_v\actson Y_v$ has a decomposition into a graph of actions of the right kind.
The theorem will follow.  

If $Y_v$ is a line, then $G_v\actson Y_v$ is either simplicial or of axial type,
so we can assume that $Y_v$ is not a line.

By Lemma \ref{lem_decompo_minimal}, 
we can first decompose $Y_v$ into a graph of actions where each 
vertex action has a global fixed point or a dense minimal subtree.
Therefore, we can assume without loss of generality that for each $v$,
$\ol G_v$ has either a global fixed point or a dense minimal subtree in $Y_v$.

\begin{lem}
  If $\ol G_v$ fixes a point in $Y_v$, then $Y_v$ is a simplicial tree.
\end{lem}

\begin{proof}
Using Lemma \ref{lem_support}, consider a finite tree $K$ spanning $Y_v$.
Let $x_0\in K$ be a $\ol G_v$-invariant point.
Denote by $K_1,\dots,K_n$ the closure
of the connected components of $K\setminus \{x_0\}$.
If $K_i,K_j$ are such that there exists $g\in \ol G_v$ with 
$g.K_i\cap K_j$ non-degenerate, replace $K_i$ and $K_j$ by $g.K_i\cup K_j$.
This way, one can assume that for all $i\neq j$ and all $g\in \ol G_v$,
$g.K_i\cap K_j=\{x_0\}$.
Since arc stabilizers are trivial, for all $g\in \ol G_v\setminus\{1\}$,
$g.K_i\cap K_i=\{x_0\}$.
Thus, $Y_v$ is the union of translates of $K_1,\dots,K_n$, all glued along $\{x_0\}$.
Since each $K_i$ is a finite tree, $Y_v$ is a simplicial tree.
\end{proof}

Now we consider the case where  $\ol G_v\actson Y_v$ has a dense minimal subtree.
Since $Y_v$ is not a line, $Y_v$ contains an infinite tripod. This tripod is fixed by $N_v$.
Moreover, consider $\ol g,\ol h\in \ol G_v$ two hyperbolic elements having distinct axes.
In particular $[\ol g,\ol h]\neq 1$. Since arc stabilizers are trivial, the axes of $\ol g$ and $\ol h$
have compact intersection. Therefore, $\ol g$ and $\ol g'=\ol g^{\ol h^k}$ have disjoint axes for $k$ large enough.
Let $g,g'$ be some preimages of $\ol g,\ol g'$ in $G$.
Then $\langle g,g'\rangle$ is a free group generated by two elements having disjoint axes in $T$, and which normalizes $N_v$.
  
\begin{lem}
  If $Y_v$ is not minimal, then $G$ splits over $N_v$.
\end{lem}

\begin{proof}
Let $G\actson S$ be the skeleton $S$ of the transverse covering $\caly$.
Since $G\actson T$ is minimal, $G\actson S$ is minimal (Lemma \ref{lem_supp_fini}).
Therefore, we need only to prove that there is an edge of $S$ whose stabilizer is $N_v$.

Let $x\in Y_v\setminus \min(Y_v)$.
Since $\min(Y_v)$ is dense in $Y_v$ and does not contain $x$,
$Y_v\setminus x$ is connected.
If $x$ does not lie in any other tree $Y\in\caly$, then $T\setminus \{x\}$ is convex,
so $T\setminus G.x$ is a $G$-invariant subtree, contradicting the minimality of $T$.

Therefore, $x$ is a vertex of $S$ and $(x,Y_v)$ is an edge of $S$ (see section \ref{sec_goa}). 
Denote by $G(x,Y_v)=\{g\in G_v| g.x=x\}$ its stabilizer.
Clearly, $N_v\subset G(x,Y_v)$.
Conversely, consider $g\in G(x,Y_v)$ and $y\in Y_v\setminus\{x\}$.
Since $y$ and $g.y$ both lie in the convex set $Y_v\setminus\{x\}$,
$[x,y]\cap [x,g.y]$ is a non-degenerate arc fixed by $g$.
Therefore, $g\in N_v$.
\end{proof}

There remains to analyse the case where  $Y_v$ 
is a not a line and $G_v\actson Y_v$ is minimal without global fix point. 
Consider the graph of groups $\Gamma=S/G$.
We identify $v$ with a vertex of $\Gamma$. 
By Lemma \ref{lem_relfg_rel} (or \ref{lem_relfg} if $G$ is finitely generated),
$G_v$ is finitely generated relative to the incident edge groups
together with at most one conjugate of each $H\in\calh$.
When $\calh\neq \es$, 
we make sure that for any  $H\in\calh$ having a conjugate in $G_v$, $\calh_v$ contains a conjugate of $H$.
Let $\ol \calh_v$ be the image of $\calh_v$ in $\ol G_v=G_v/N_v$.
Then we can apply the relative version of Sela's Theorem (Proposition \ref{prop_sela_rel})
to $(\ol G_v,\ol\calh_v)\actson Y_v$.

Assume first that $(\ol G_v,\ol\calh_v)$ is freely decomposable.
Therefore, $G_v$ splits over $N_v$ relative to $\calh_v$.
Thus, we can refine $\Gamma$ at $v$ using this splitting,
so $G$  splits over $N_v$. 
This is really a splitting of $(G,\calh)$ because 
we made sure that for each $H\in\calh$,
either $H$ has a conjugate in $\calh_v$, or $H$ is conjugate in some other vertex group of $\Gamma$.

Assume now that $\ol G_v\actson Y_v$ has a decomposition into
a graph of actions where each vertex action is either
axial or of surface type.
Then clearly, the action $G_v\actson Y_v$ as a decomposition into a graph of actions
where each vertex action is either axial or of Seifert type.
This completes the proof of Theorem \ref{thm_pw2triv}.
\end{proof}

\section{Proof of Main Theorem and corollaries}\label{sec_proof}


Recall that an action of the pair $(G,\calh)$ on a tree is an action of $G$ where
each $H\in\calh$ fixes a point. In terms of graphs of groups, a splitting of the pair $(G,\calh)$
is an isomorphism of $G$ with a graph of groups such that each $H_i$
is contained in a conjugate of a vertex group (see section \ref{sec_relfg}).

\begin{thm}[Main Theorem]\label{thm_main}
Consider a minimal action of finitely generated pair $(G,\calh)$ on an $\bbR$-tree $T$ by isometries.
Assume that
\begin{enumerate*}
\item $T$ satisfies the ascending chain condition; 
\item there exists a finite family of arcs $I_1,\dots,I_p$ such that $I_1\cup\dots\cup I_p$ spans $T$ (see Definition \ref{dfn_span})
and such that for any unstable arc $J$ contained in some $I_i$,
\begin{enumerate*}
\item $G(J)$ is finitely generated;
\item $G(J)$ is not a proper subgroup of any conjugate of itself
\ie $\forall g\in G$, $G(J)^g\subset G(J)\Rightarrow G(J)^g= G(J)$.
\end{enumerate*}
\end{enumerate*}

Then either $(G,\calh)$ splits over the stabilizer of an unstable arc, or over the stabilizer of an infinite tripod 
(whose normalizer contains a non-abelian free group generated by two elements having disjoint axes),
or $T$ has a decomposition into a graph of actions where each
vertex action is either 
\begin{enumerate*}
\item simplicial: $G_v\actson Y_v$ is a simplicial action;
\item of Seifert type: the vertex action $G_v\actson Y_v$ has kernel $N_v$, and the faithful action $G_v/N_v\actson Y_v$
is dual to an arational measured foliation on a closed $2$-orbifold with boundary;
\item axial: $Y_v$ is a line, and the image of $G_v$ in $\Isom(Y_v)$ is a finitely generated group 
acting with dense orbits on $Y_v$.
\end{enumerate*}
\end{thm}

\begin{rem*}
The group over which $G$ splits (\ie the stabilizer of an unstable arc or of a tripod)
is really its full pointwise stabilizer. This contrasts with \cite[Theorem 9.5]{BF_stable}.
\end{rem*}

\begin{proof}[Proof of Main Theorem]
  Let $(G,\calh)\actson T$ be as in Main Theorem. 
By Theorem \ref{thm_acc2pw}, if $(G,\calh)$ does not split over the stabilizer of an unstable arc,
$T$ is piecewise stable.
By Theorem \ref{thm_pw2triv}, either $(G,\calh)$ splits over the stabilizer of an infinite tripod with the required properties, 
or $T$ has a decomposition into a graph of actions of the desired type.
\end{proof}

\begin{SauveCompteurs}{corBF}
  \begin{cor}\label{cor_BF}
  Under the hypotheses of Main Theorem,
either $T$ is a line or $(G,\calh)$ splits over a subgroup $H$ which is an extension
of a cyclic group (maybe finite or trivial) by an arc stabilizer.
\end{cor}
\end{SauveCompteurs}

\begin{proof}
We can assume that $T$ splits as a graph of actions $\calg$ where each vertex action $G_v\actson Y_v$ is either simplicial,
of Seifert type, or axial.
Let $(Y_v)_{v\in V}$ be the family of non-degenerate vertex trees.
Let $S$ be the skeleton of this transverse covering (see Lemma \ref{lem_transverse_cov}).

First, consider the case where $S$ is reduced to a point $v$.
This means that $T=Y_v$, and $G=G_v$, and that 
$T$ is itself simplicial or of Seifert type ($T$ is not a line so cannot be of axial type).
If $T$ is simplicial, the result is clear, so assume that $T$ is of Seifert type.
Let $N_v\normal G_v$ be the kernel of this action and let $\Sigma$ be a $2$-orbifold with boundary, with cone singularities,
 such that $G_v/N_v=\pi_1(\Sigma)$ and holding an arational measured foliation to which $G_v/N_v\actson Y_v$ is dual.
Consider a splitting of $G_v/N_v$ corresponding to an essential simple closed curve 
(\ie a curve which cannot be homotoped to a point, a cone point, or to the boundary).
Such a curve exists because $\Sigma$ holds an arational measured foliation.
This defines a splitting of $\pi_1(\Sigma)$ over a cyclic group, and a splitting of $G=G_v$ 
 over an extension of a cyclic group by $N_v$.
Any subgroup $H\subset G$ elliptic in $T$ corresponds to the fundamental group of
a cone point or of a boundary component of $\Sigma$.
Thus $H$ is elliptic in this splitting and we get a splitting of the pair $(G,\calh)$.


Now assume that $S$ is not reduced to a point.
By Lemma \ref{lem_supp_fini}, $G\actson S$ is minimal and any $H\in\calh$ is elliptic in $S$.
In particular, given any edge $e$ of $S$, the corresponding splitting of $G$ over $G(e)$ is non-trivial.
We shall prove that for some edge $e$ of $S$, $G(e)$ is an extension of an arc stabilizer by a cyclic group,
and the corollary will follow.

Assume that some action $G_v\actson Y_v\simeq\bbR$ is axial.
Consider an edge $e=(x,Y_v)$ of $S$ incident on $v$ (see Lemma \ref{lem_transverse_cov}).
Its stabilizer $G(e)$ is the stabilizer of $x$ in $G_v$.
Since the stabilizer of $x$ in $G_v/N_v$ is either trivial or $\bbZ/2\bbZ$,
$G(e)$ is a extension of a cyclic group by $N_v$. 

Assume that some action $G_v\actson Y_v$ is of Seifert type. 
Consider an edge $e=(x,Y_v)$ incident on $v$ in $S$.
Since $G(e)$ is the stabilizer of $x$ in $G_v$,
$G(e)$ is an extension by $N_V$ of the stabilizer of a point in $\pi_1(\Sigma)\actson Y_v$ which is cyclic.

The only remaining case is when each $G_v\actson Y_v$ is simplicial. In this case $G\actson T$ is simplicial,
and the result is clear.
\end{proof}

\begin{SauveCompteurs}{cor_petit}
\begin{cor}\label{cor_small}
  Let $G$ be a finitely generated group for which any small subgroup is finitely generated.
Assume that $G$ acts on an $\bbR$-tree $T$ with small arc stabilizers. 

Then $G$ splits over the stabilizer of an unstable arc, or over a tripod stabilizer,
or $T$ has a decomposition into a graph of actions as in Main Theorem.
In particular, $G$ splits over a small subgroup.
\end{cor}
\end{SauveCompteurs}

\begin{proof}
  We prove that all hypotheses of the Main Theorem are satisfied.
The set of small subgroups is closed under increasing union.
Since small subgroups are finitely generated, any ascending chain of small subgroups is finite.
The ascending chain condition follows. Similarly, if $H\subsetneq H^g$ for some small subgroup $H\subset G$,
then the set of small subgroup $H^{g^n}$ is an infinite increasing chain, a contradiction.
Then Main Theorem and corollary \ref{cor_BF} apply.
\end{proof}

In the following situation, stabilizers of unstable arcs can be controlled by tripod stabilizers.
\begin{SauveCompteurs}{cor_tripod}
\begin{cor}\label{cor_tripod}
Consider a finitely generated group $G$ acting by isometries on an $\bbR$-tree $T$.
Assume that
\begin{enumerate*}
\item arc stabilizers have a nilpotent subgroup of bounded index (maybe not finitely generated);
\item tripod stabilizers are finitely generated  (and virtually nilpotent);
\item no group fixing a tripod is a proper subgroup of any conjugate of itself;
\item any chain $H_1\subset H_2\dots$ of tripods stabilizers stabilizes. 
\end{enumerate*}

Then  either $G$ splits over a subgroup having a finite index subgroup fixing a tripod, 
or $T$ has a decomposition as in the conclusion of Main Theorem.\\
\end{cor}
\end{SauveCompteurs}

\begin{proof}
Let $k$ be a bound on the index of a nilpotent subgroup in an arc stabilizer.
We first claim that the stabilizer of an unstable arc has a subgroup of bounded index
fixing a tripod.
Indeed, let $A$ be the stabilizer of an unstable arc $I$, and let $B\supsetneq A$ be the stabilizer of a sub-arc $J$.
If $A$ has index at most $k$ in $B$, then $A$ contains a normal subgroup $A'$ of index at most $k!$.
Since for $g\in B\setminus A$, $g.I$ is an arc distinct from $I$ and containing $J$,
$B.I$ contains a tripod. This tripod is fixed by $A'$.
Suppose that the index of $A$ in $B$ is larger than $k$. 
Let $N_B$ be a nilpotent subgroup of index at most $k$ in $B$ and $N_A=A\cap N_B$.
Because of indices, $N_A\subsetneq N_B$. Let $N$ be the normalizer of $N_A$ in $N_B$.
It is an easy exercise to check that in a nilpotent group, no proper subgroup is its own normalizer
so $N_A\subsetneq N$.
Therefore, $N.I$ contains a tripod, and this tripod is fixed by $N_A$.

It follows that the stabilizer of an unstable arc is finitely generated.
Since each tripod stabilizer is slender, any ascending chain of subgroups fixing tripods stabilizes.
It follows that the stabilizer of an unstable arc cannot be properly conjugated into itself.
Consider $A_1\subsetneq\dots\subsetneq A_n\subsetneq$ an ascending chain of arc stabilizers.
Let $N_n$ be a subgroup of bounded index of $A_n$ fixing a tripod, and let $N_{i,n}$ be its intersection with $A_i$
for $i\leq n$. Since $N_{i,n}$ has bounded index in $A_i$, it takes finitely many values for each $i$. By a diagonal argument,
we get an ascending chain of subgroups fixing tripods, a contradiction.
Thus, Main Theorem applies, and the corollary is proved.
\end{proof}

\section{An example}\label{sec_example}

The goal of this section is to give an example of an action with trivial tripod stabilizers
providing a counter-example to \cite[Theorem 10.8]{RiSe_structure} and 
\cite[Theorem 2.3 and 3.1]{Sela_acylindrical}.
Note however that the proof of Theorem 3.1 in \cite{Sela_acylindrical}
is valid under the stronger assumption that $T$ is \emph{super-stable},
which is a natural hypothesis since it is satisfied in applications in \cite{RiSe_structure, Sela_acylindrical,Sela_diophantine1}.

\begin{thm}
  There exists a non-trivial minimal action of a group $G$ on an $\bbR$-tree $T$ such that
  \begin{itemize*}
  \item $G$ is finitely presented (and even word hyperbolic) and freely indecomposable
  \item tripod stabilizers are trivial
  \item $T$ satisfies the ascending chain condition
  \item $T$ has no decomposition into a graph of actions as in 
Main Theorem.
\end{itemize*}
\end{thm}

\begin{proof}
Let $A$ be a freely indecomposable hyperbolic group containing a free malnormal
subgroup $M_1=\langle a,b_1\rangle$. For instance, one can take for $A$ a surface group,
and for $M_1$, the fundamental group of a punctured torus contained in this surface.
Let $A'$ be another copy of the group $A$, $C=\langle a\rangle$, and let $G=A*_C A'$.
By Bestvina-Feighn's Combination Theorem, $G$ is hyperbolic \cite{BF_combination}.

For $i>1$, define inductively $M_i\subset M_{i-1}$ by
$M_i=\langle a,b_i\rangle$ where $b_i=b_{i-1}ab_{i-1}^2$.
One can easily check that $M_{i}$ is malnormal in $M_{i-1}$ (one can also use the
software Magnus to do so \cite{Magnus_software}).
Moreover, one also easily checks that $\cap_{i\geq 1} M_i=C$
since any reduced word $w$ on $\{a,b_i\}$ 
defines a reduced word on $\{a,b_{i-1}\}$ by the obvious substitution, and
its length is strictly larger if $w$ is not a power of $a$.

Let $\Gamma_k$ be the graph of groups 
and $G\actson T_k$ be the corresponding Bass-Serre tree as shown in Figure \ref{fig_exple}.

\begin{figure}[htb]
  \centering
  \includegraphics{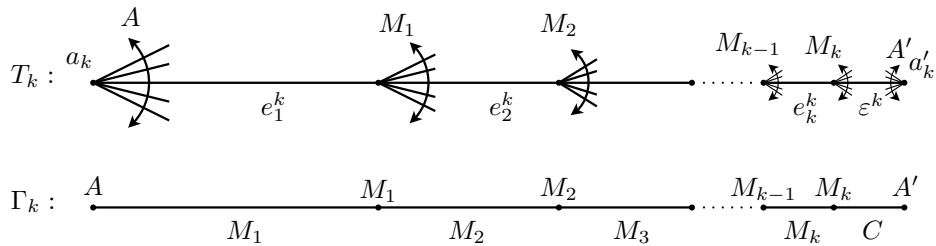}
    \caption{The graph of groups $\Gamma_k$ and the tree $T_k$\label{fig_exple}}
\end{figure}

Let $a_k,a'_k$ be the vertices of $T_k$ fixed by $A$ and $A'$ respectively.
Denote by $e^k_i\subset[a_k,a'_k]$ the edge of $T_k$ such that $G(e^k_i)=M_i$, and by $\eps^k\subset [a_k,a'_k]$ the edge such that $G(\eps^k)=C$.
Figure \ref{fig_exple} shows a neighbourhood of $[a_k,a_k']$ in $T_k$, and the action of the vertex stabilizers.
The main feature is that from each vertex $v$ with stabilizer $M_i$, the neighbourhood of $v$ consists of one orbit
of edges $M_i.e_{i+1}^k$, and one single edge $e_i^k$.

Assign length $1/2^i$ to the edge $e_i^k$, and length $1/2^k$ to $\eps^k$.
This way, $d(a_k,a'_k)=1$ for all $k$.
There is a natural morphism of $\bbR$-trees $f_k:T_k\ra T_{k+1}$ sending $e^k_i$ to $e^{k+1}_i$ and 
$\eps^k$ to $e_{k+1}^{k+1}\cup\eps^{k+1}$.


The length function of $T_{k+1}$ is not larger
than the length function of $T_{k}$, so $T_k$ converges in the length function topology to an action on an $\bbR$-tree $T$.
This action is non-trivial since for any $g\in A\setminus C$ and $g'\in A'\setminus C$, 
the translation length of $gg'$ in every $T_k$ is $2$.

We shall prove that $T_k$ converges strongly to $T$.
Consider $\alpha,\beta$ two distinct edges of $T_0$ sharing a vertex $v$
($T_0$ is the tree dual to the amalgam $G=A*_C A'$).
If $v$ is in the orbit of $a'_0$, then the image of $\alpha$ and $\beta$ in $T_k$
share only one point. Now assume without loss of generality that $v=a$, and $\alpha=\eps^0$.
Let $g\in A\setminus C$ be such that $g.\alpha=\beta$.
Consider $k_0$ the smallest integer such that, $g\notin M_{k_0}$. Then for $k\geq k_0$,
the image of $\alpha$ and $\beta$ in $T_k$ share exactly an initial segment of length
$1-1/2^{k_0}$. In particular, $f_k$ is one-to-one in restriction to the union of the images
of $\alpha$ and $\beta$.
Now if $K$ is a finite subtree of $T_0$, it follows that for $k$ large enough, $f_k$
is one-to-one in restriction to the image of $K$ in $T_k$.
This proves the strong convergence of $T_k$.

Say that an arc $I\subset T_k$ is \emph{immersed in $\Gamma_k$}
if the restriction to $I$ of the quotient map $T_k\ra \Gamma_k$ is an immersion.
We claim that for any arc $I\subset T_k$ which is not immersed in $\Gamma_k$,
$G(I)$ is trivial. It is sufficient to prove it for $I$ of the form
$I=\alpha\cup \beta$ where $\alpha,\beta\in E(T_k)$ are incident on a common vertex
and are in the same orbit. The triviality of $G(I)=G(\alpha)\cap G(\beta)$
follows easily from the malnormality of $M_1$ in $A$, of $M_{i+1}$ in $M_i$,
and of $C$ in $M_k$ and $A'$.
Since a tripod cannot be immersed, it follows that tripod stabilizers of $T_k$ are trivial.
Going to the limit, tripod stabilizers of $T$ are trivial.

Now we study arc stabilizers of $T$.
Denote by $F_k:T_k\ra T $ the natural map, and by $a=F_k(a_k)$ and $a'=F_k(a'_k)$.
Let $I=[u,v]\subset T$ be a non-degenerate arc with non-trivial stabilizer. 
By strong convergence, if $g\in G(I)$ then there exists a lift $I_k\subset T_k$
in restriction to which $F_k$ is isometric, and which is fixed by $g$.
The argument above shows that the image of $I_k$ in $\Gamma_k$ is immersed.
In particular, we can assume that $I_k\subset [a_k,a'_k]$. 
If $I$ contains $a'$, then since $F_k\m(a')=\{a'_k\}$, $I_k$ intersects $\eps_k$ in a non-degenerate segment,
so $G(I_k)=C$ and $G(I)=C$.
If $I$ does not contain $a'$, for $k$ large enough, $I_k$ does not intersect $\eps_k$
so $G(I_k)=M_i$ for some $i$ independent of $k$.
If follows that $G(I)=M_i$.
The ascending chain condition for $T$ follows immediately.

Since the restriction of $F_k$ to the ball of radius $1-1/2^k$ around $a_k$ 
is an isometry, the segment $[a,a']\subset T$ contains an infinite number of branch points.
In particular, $T$ is not a simplicial tree.

We claim that no subtree $Y$ of $T$ is indecomposable.
Since $[a,a']$ spans $T$, we may assume that $Y\cap [a,a']$ contains a non-degenerate arc $I$.
By indecomposability, the orbit of any point of $I$ is dense in $I$, a contradiction.

Finally, assume that $T$ has a decomposition as in Main Theorem. 
Since there are no indecomposable subtrees in $T$, $T$ is simplicial, a contradiction.
\end{proof}

\appendix
\section{Sela's Theorem}\label{sec_Sela}

\begin{SauveCompteurs}{thm_sela}
\begin{thm}[\cite{Sela_acylindrical}]\label{thm_Sela}
  Consider a minimal action of a finitely generated group $G$ on an $\bbR$-tree $T$ with trivial arc stabilizers.

Then, either $G$ is freely decomposable, or $T$ has a decomposition into a graph of actions where each
non-degenerate vertex action is either 
\begin{enumerate*}
\item of surface type: the vertex action $G_v\actson Y_v$ 
is dual to an arational measured foliation on a closed $2$-orbifold with boundary;
\item axial: $Y_v$ is a line, and the image of $G_v$ in $\Isom(Y_v)$ is a finitely generated group 
acting with dense orbits on $Y_v$.
\end{enumerate*}
\end{thm}
  \end{SauveCompteurs}


\subsection{Preliminaries}

The following lemma is essentially contained in \cite{LP}.
\begin{lem}\label{lem_cv2}
Let $G\actson T$ be a minimal action with trivial arc stabilizers.

There exists a sequence of $\bbR$-trees $G_k\actson T_k$
converging strongly to $G\actson T$ such that
$T_k$ is geometric, dual to a foliated 2-complex $G_k\actson X_k$ whose leaf space is Hausdorff,
and such that $\Phi_k:G_k\ra G$ is one-to-one in restriction to each point stabilizer.
\end{lem}

\begin{proof}
  Consider a sequence of actions $G_k\actson T_k$ 
dual to foliated $2$-complexes $X_k$ as in Theorem \ref{thm_strongCV}.
Let $N_k$ be the kernel of $\Phi_k:G_k\ra G$.
Let $L_k$ be the subgroup of $N_k$ generated by elements which preserve a leaf of $X_k$.
This is a normal subgroup of $G_k$.
Let $X'_k=X_k/L_k$ be the quotient foliated $2$-complex, endowed with the natural action of $G'_k=G_k/L_k$.
By construction, for each leaf $l$ of $X'_k$, the global stabilizer of $l$ embeds into $G$
under the induced morphism $\Phi'_k:G'_k\ra G$.

The natural map $X_k\ra X'_k$ is a covering.
Since $\pi_1(X_k)$ is generated by free homotopy classes of curves contained in leaves,
and since $L_k$ is generated by elements preserving a leaf,
$\pi_1(X'_k)$ generated by free homotopy classes of curves contained in leaves.
Let $T'_k$ be the leaf space made Hausdorff of $X'_k$. It is an $\bbR$-tree 
by Theorem \ref{thm_LP_leaves}.
Since $L_k\subset N_k$, the map $f_k:T_k\ra T$ factors through the natural map
$T_k\ra T'_k$. It follows that $T'_k$ is geometric, dual to $X'_k$, and that $G'_k\actson T'_k$ converge
strongly to $G\actson T$.

By Proposition \ref{prop_LP_separation},
there exists a countable union of leaves $\cals$ on the complement of which
two points of $X'_k$ are identified in $T'_k$ if and only if they lie on the same leaf.
It follows that any element $g\in G'_k$ fixing an arc in $T'_k$ preserves a leaf of $X'_k$.
Since arc stabilizers of $T$ are trivial, the image of $g$ in $G$ is trivial, so $g=1$.
Thus, arc stabilizers of $G'_k\actson T'_k$ are trivial.

By proposition \ref{prop_hausdorff} below, the leaf space of $X'_k$ is Hausdorff. 
In particular, a point stabilizer of $T'_k$ coincides with a leaf stabilizer of $X'_k$,
which embeds into $G$.
The lemma follows.
\end{proof}

\begin{dfn}
We will say that $G\actson T$ is \emph{nice} if
\begin{itemize*}
\item $T$ is geometric, dual to a foliated $2$-complex $G\actson X$ such that
$\pi_1(X)$ is generated by free homotopy classes of curves contained in leaves;
\item arc stabilizers of $T$ are trivial.
\end{itemize*}
\end{dfn}

The following result follows from the concatenation of Lemma 3.5 and 3.4 in \cite{LP}:
\begin{prop}[\cite{LP}]\label{prop_hausdorff}
  Assume that $T$ is nice.

Then the leaf space of $X$ is Hausdorff: any two points of $X$ identified in $T$
are in the same leaf.
\end{prop}

\begin{lem}
  Assume that $G\actson T$ is nice.
Let $Y$ be an indecomposable component of $T$ (as in Proposition \ref{prop_decompo}), 
and $H$ its global stabilizer.

Then $H\actson Y$ is nice.
\end{lem}

\begin{proof}
Let $G\actson X$ be a foliated $2$-complex such that $T_k$ is dual to $X$.
Let $C\subset X/G$ be its cut locus (see Definition \ref{dfn_cut}), and $\Tilde C$ its preimage in $X$.
The tree $Y$ is dual to the closure $\ol U$ of a connected component $U$ of $X\setminus \Tilde C$.
By Remark \ref{rem_decompo_geom},  $\pi_1(\ol U)$
is generated by free homotopy classes of curves contained in leaves.
This means that $Y$ is nice.  
\end{proof}

\begin{prop}\label{prop_3types}
Assume that $G\actson T$ is nice and let $H\actson Y$ be an indecomposable component.
Then one of the following holds:
\begin{itemize*}
\item Axial type: $Y$ is a line, and the image of $H$ in $\Isom(Y)$ is a finitely generated group 
acting with dense orbits on $Y$;
\item Surface type: $H$ is the fundamental group of a $2$-orbifold with boundary $\Sigma$ holding an arational measured foliation 
and $Y$ is dual to $\Tilde \Sigma$;
\item Exotic type: $H$ has a non-trivial free decomposition $H\actson S$ in which any subgroup of $H$ fixing a point in $Y$ 
fixes a point in $S$.
\end{itemize*}
\end{prop}

\begin{proof}
This proposition is essentially well known: it is a way of describing the output of the Rips machine.
This would follow from \cite{BF_stable}
if we knew that $H$ is finitely presented (because of the finiteness hypothesis in \cite[Definition 5.1]{BF_stable}).
But we can apply some arguments of \cite{Gui_approximation} where
finite presentation is not assumed.

We recall some vocabulary from \cite{GLP1} or \cite{Gui_approximation}.
A closed multi-interval $D$ is a finite union of compact intervals.
A \emph{partial isometry} of $D$ is an isometry between closed sub-intervals of $D$.
A \emph{system of isometries} $\cals$ on $D$
is a finite set of partial isometries of $D$.
Its \emph{suspension} is the foliated $2$-complex obtained from $D$
by gluing for each partial isometry $\phi:I\ra J\in\cals$ by 
a foliated band $I\times [0,1]$ on $D$ joining $I$ to $J$ whose holonomy is given by $\phi$.
A \emph{singleton} is a partial isometry $\phi:I\ra J$ where $I$ is reduced to a point.
We also call singleton the band corresponding to a singleton in the suspension of a system of isometries.
The $\rond\cals$-orbits are the equivalence classes for the equivalence relation generated by
$x\sim y$ if $y=\phi(x)$ for some non-singleton $\phi:I\ra J$ such that $x\in \rond I$.
A system of isometries $\cals$ is \emph{minimal} if $\cals$ has no singleton and
every $\rond\cals$-orbit of any point in $D\setminus \partial D$ is dense.
The system $\cals$ is \emph{simplicial} if every $\cals$-orbit is finite.

The proposition holds if $H\actson Y$ is dual to a minimal system of isometries: 
this follows from \cite{Gui_approximation} section 5 (axial case), proposition 7.2 (exotic case), 
and section 8 (surface case).

Now $H\actson Y$ is dual to some foliated $2$-complex $H\actson X$, occurring as
a Galois covering of the suspension of some system of isometries $\cals$.
Following \cite{GLP1}, \cite{BF_stable} or \cite[Prop 3.1]{Gui_approximation},
one can perform a sequence of Rips moves on $\cals$
so that $\cals$ becomes a disjoint union finitely many systems of isometries 
which are either minimal or simplicial, together with a finite set of singletons joining them.
This decomposition of induces a decomposition of $Y$ as a graph of actions as in Proposition \ref{prop_decompo}.
Since $Y$ is indecomposable, this decomposition is trivial, and $Y$ is dual to the suspension of a minimal component of $\cals$.
\end{proof}

\subsection{Proof of the theorem}

\begin{proof}[Proof of Sela's Theorem]

Let $G_k\actson T_k$ be a sequence of geometric actions converging strongly to $G\actson T$
as in lemma \ref{lem_cv2}. Arc stabilizers of $G_k\actson T_k$ are trivial.
By Proposition \ref{prop_decompo}, $T_k$ splits as a graph of actions $\calg$, where each non-degenerate vertex action
is either simplicial or indecomposable. 

Assume that for all $k$, the simplicial part of $T_k$ is non-empty. We shall prove that $G$ is freely decomposable.
Collapse the indecomposable components of $T_k$ as in Definition \ref{dfn_collapse}.
Let $S_k$ be the obtained tree. Clearly, $S_k$ is a simplicial tree, and edge stabilizers are trivial
because arc stabilizer of $T_k$ are trivial.
The action $G_k\actson S_k$ is minimal by Lemma \ref{lem_collapse} and therefore non-trivial.
The map $f_k:T_k\ra T_{k+1}$ maps an indecomposable component of $T_k$ into an indecomposable component of $T_{k+1}$
so $f_k$ induces a map $S_k\ra S_{k+1}$.
By  Scott's Lemma, $G$ is freely decomposable.

Now, we assume that for all $k$, $T_k$ splits as a graph of indecomposable actions.
Let $S_k$ be the skeleton of the corresponding transverse covering of $T_k$.
Recall that its vertex set $V(S)$ is $V_0(S)\cup V_1(S)$ where $V_1(S)=\caly$, and $V_0(S)$ is the set of
points $x\in T$ lying in the intersection of two distinct trees of $\caly$.
Since $f_k$ maps an indecomposable tree into an indecomposable tree, $f_k$ induces a map $V_1(S_k)\ra V_1(S_{k+1})$.
Moreover, $f_k$ induces a map $V_0(S_k)\ra V_0(S_{k+1})\cup V_1(S_{k+1})$:
for $x\in V_0(S_k)$, if $f_k(x)$ belongs to two distinct indecomposable components,
we map $x$ to $f_k(x)\in V_0(S_{k+1})$,
otherwise, we map $x$ to the only indecomposable component containing $f_k(x)$.
This map extends to a map $g_k:S_k\ra S_{k+1}$ sending an edge to an edge or a vertex.
The number of orbits of edges of $S_k$ is non-increasing so
for $k$ large enough, no edge of $S_k$ is collapsed by $g_k$,
and any pair of edges $e_1,e_2$ identified by $g_k(e_1)$ are in the same orbit.
Moreover, using Scott's Lemma, we can assume that for $k$ large enough, no edge of $S_k$ has trivial stabilizer.

The following lemma will be proved in next sections.

\begin{lem}[Stabilisation of indecomposable components]\label{lem_inj}
  For $k$ large enough, the following holds.
Consider an indecomposable component $H\actson Y$ of $T_k$, and let $H'\actson Y'$ be the 
indecomposable component of $T_{k+1}$ containing $f_k(Y)$.

Then $f_k{}_{|Y}$ is an isometry from $Y$ onto $Y'$.
Moreover, if $Y$ is not a line, then $\phi_k(H)=H'$.
\end{lem}

Using this lemma, we shall prove that for $k$ large enough, $f_k:T_k\ra T_{k+1}$ is an isometry.
If this is not the case, then there exists two arcs $J_1,J_2\subset T_k$
with $J_1\cap J_2=\{x\}$ and $f_k(J_1)=f_k(J_2)$.
By shortening them, we may assume that $J_1$ and $J_2$ lie in some indecomposable
components $Y_1$, $Y_2$. By Lemma \ref{lem_inj}, $Y_1\neq Y_2$.
Therefore, the edges $(x,Y_1)$ and $(x,Y_2)$ of $S_k$ are identified under $g_k$. 
By the assumption above, they lie in the same orbit.
Consider $g\in G_k(x)$ such that $g.Y_1=Y_2$.
Let $Y'=g_k(Y_1)=g_k(Y_2)$. 

First, assume that $Y_1$ is not a line.
Since $\phi_k(g)$ preserves $Y'$, 
there exists $\Tilde g\in G_k(Y_1)$ with $\phi_k(\Tilde g)=\phi_k(g)$ by Lemma \ref{lem_inj}.
Since $f_k$ is isometric in restriction to $Y_1$ and since $\phi_k(\Tilde g)=\phi_k(g)$ fixes $f_k(x)$, 
$\Tilde g$ fixes $x$.
In particular, $g\m \Tilde g$ fixes $x$ and lies in the kernel of $\phi_k$.
By Lemma \ref{lem_cv2}, $\phi_k$ is one-to-one in restriction to point stabilizers, so $g=\Tilde g$. 
This is a contradiction because $g \notin G_k(Y_1)$.

Therefore, $Y_1$ is a line.
Since $\phi_k$ is one-to-one in restriction to $G_k(x)$,
$\phi_k(g)\neq 1$. 
Since $\phi_k(g)$ preserves the line $Y'$, and since arc stabilizers are trivial, 
$\phi_k(g)$ acts on $Y'$ as the reflection fixing $f_k(x)$. 
By the assumption above, edge stabilizers of $S_{k}$ are non-trivial
so consider $h\neq 1$ in the stabilizer of the edge $(x,Y_1)$.
Since $h$ fixes a point in $T_k$ (namely, $x$), $\phi_k(h)\neq 1$.
So $\phi_k(h)$ acts as the same reflection as $\phi_k(g)$ on $Y'$.
Since arc stabilizers are trivial, $\phi_k(g)=\phi_k(h)$.
Since $g,h$ are contained in a point stabilizer, we get $g=h$. 
This is a contradiction since $g\notin G_k(Y_1)$.

This proves that for $k$ large enough, $f_k$ is an isometry.
It follows that $\phi_k$ is an isomorphism because any element
of $\ker \phi_k$ must fix $T_k$ pointwise.

Thus, Proposition \ref{prop_3types} gives us a decomposition of $T$.
If there is an exotic type vertex action $H\actson Y$,
the free splitting of $H$ given in Proposition \ref{prop_3types} 
can be used to refine the decomposition of $G$ induced by the skeleton $S_k$, 
thus giving a decomposition of $G$ as a free product.

Finally, there remains prove the finite generation claimed in the axial case.
Let $H\actson Y$ be an indecomposable component such that $Y$ is a line.
In particular, $H$ occurs as a vertex group in the graph of groups given by the decomposition
of $T$. 
By Lemma \ref{lem_relfg}, $H$ is finitely generated relative to the stabilizers of incident edges
 $H_1,\dots, H_p$. 
Denote by $\psi:H\ra \Isom(Y)$ the map induced by the action.
Since each $H_i$ fixes a point in $Y$, $\psi( H_i)$ is cyclic of order at most 2.
Since $\psi(H)$ is finitely generated relative to $\psi(H_1),\dots, \psi( H_p)$,
$\psi(H)$ is finitely generated.
\end{proof}

\subsection{Standard form for indecomposable components}

We set up the material needed for the proof of Lemma \ref{lem_inj}.
Consider $H\actson Y$ a nice and indecomposable action.
Given an arc $I\subset Y$ and a finite generating set $S\subset H$, 
we construct a foliated $2$-complex $X(I,S)$ as a kind of thickened Cayley graph.
Alternatively, one can view $X(I,S)$ as the covering space with deck group $H$ of
 the suspension of a system of isometries on $I$.
The point here is to start with an arbitrary arc $I$ and to allow $S$ to be large in order to ensure that $Y$ is dual 
to $X(I,S)$.

We first define $X(I,S)$ for any arc $I\subset Y$, and any finite subset $S\subset H$.
Start with $H\times I$, endowed with the action of $H$ given by $g.(h,x)=(gh,x)$.
For each $s\in S$, consider $K_s=(I\cap s\m I)$, and
for each $g\in H$, add a foliated band $K_s\times [0,1]$
joining $\{g\}\times K_s$ to $\{gs\}\times s.K_s$
whose holonomy is given by the restriction of $s$ to $K_s$.
The map sending $(g,x)$ to $g.x$ extends uniquely to a map $p:X(I,S)\ra Y$
which is constant on each leaf.

\begin{lem}\label{lem_dual}
Consider $H\actson Y$ a nice indecomposable action.
  For every arc $I\subset Y$, there exists a finite generating set $S$ of $H$
such that 
\renewcommand{\theenumi}{\roman{enumi}}
\renewcommand{\labelenumi}{\textup{(\theenumi)}}
\begin{enumerate*}
\item the leaf space of $X(I,S)$ is Hausdorff;
\item $Y$ is dual to $X(I,S)$.
\end{enumerate*}
\end{lem}
 
Before proving the lemma, we introduce some tools.
A \emph{holonomy band} 
in a foliated $2$-complex $\Sigma$ is a continuous map $b:I\times [0,1]\ra \Sigma$
such that $b_{|I\times \{0\}}$ (resp. $b_{|I\times \{1\}}$) is an isometric map to a subset of a transverse edge of $\Sigma$,
and  for each $x\in I$, $b(\{x\}\times [0,1])$ is contained in a leaf segment, and this leaf segment is regular if $x\notin \partial I$.
The following property is a restatement of Theorem 2.3 in \cite{GLP2}.

\begin{prop}[Segment-closed property, \cite{GLP2}]
Let $\Sigma$ be a compact foliated $2$-complex.
Consider two arcs $I,J\subset \Sigma^{(1)}$ and an isometry $\phi:I\ra J$ 
such that for all but countably many $x\in I$,
$x$ and $\phi(x)$ are in the same leaf.

Then there exist a subdivision $I=I_1\cup \dots\cup I_p$,
and for each $i\in\{1,\dots,p\}$ a holonomy band $b_i$ joining $I_i$ to $\phi(I_i)$
whose holonomy is the restriction of $\phi$.\qed
\end{prop}

\begin{cor}\label{cor_bands}
Assume that $G\actson T$ is nice,
 dual to a foliated $2$-complex $G\actson X$.

Consider two arcs $I,J$ contained in transverse edges of $X$ and having the same image in $T$.
Then there exist subdivisions $I=I_1\cup \dots\cup I_p$, $J=J_1\cup \dots\cup J_p$,
such that for each $i\in\{1,\dots,p\}$ there exists a holonomy band joining $I_i$ to $J_i$.
\end{cor}

\begin{proof}
  Let $\phi:I\ra J$ be an isometry such that for all $x\in I$, $x$ and $\phi(x)$ map to the same point in $T$.
By Proposition \ref{prop_LP_separation}, for all but countably many $x\in I$,
$x$ lies in the same leaf as $\phi(x)$.
Apply segment closed property to the images of $I,J$ in the compact foliated $2$-complex $X/G$.
By lifting the obtained holonomy bands to $X$, we obtain
a subdivision $I=I_1\cup \dots\cup I_p$ and for each $i\in\{1,\dots,p\}$ 
a holonomy band $b_i$ joining $I_i$ to $g_i.\phi(I_i)\subset g_i.J$ for some $g_i\in G$.
Since $x$ and $\phi(x)$ map to the same point in $T$, $g_i$ fixes the arc $I_i$,
so $g_i=1$. 
\end{proof}

\begin{proof}[Proof of Lemma \ref{lem_dual}]
Consider a foliated $2$-complex $H\actson X$ such that $Y$ is dual to $X$.
Let $\pi:X\ra X/H$ be the covering map.
Consider the cut locus $C\subset X/H$,
$\Tilde C$ its preimage in $X$. The set 
$U=X\setminus \Tilde C$ is connected: 
otherwise, $Y$ would split as a graph of actions with at least two non-degenerate vertex trees
by Proposition \ref{prop_decompo}; this is impossible because $Y$ is indecomposable (Lemma \ref{lem_indec_component}).
Let $q:X\ra Y$ the quotient map ($Y$ is the leaf space of $X$).

Recall that $\Xt\subset X^{(1)}$ denotes the union of all closed transverse edges of $X$.
We aim to construct a collection of partial isometries $\Phi$
 from $\Xt$ to $H\times I$, and use $\Phi$ to transport the holonomy of the triangles of $X$ 
to bands on $H\times I$. 

Consider a finite set of arcs $J_1,\dots, J_p$ contained in edges of $\Xt$
such that $H.(J_1\cup\dots\cup J_p)=\Xt$
and $q(J_1)\cup\dots \cup q(J_p)\supset I$.
By subdividing each $J_i$, one may assume that either $q(J_i)\subset I$, or $q(J_i)\cap I$ is degenerate.
Since $Y$ is indecomposable, $I$ spans $Y$, so up to further subdivision of each $J_i$, one may assume that
for each $i\in\{1,\dots,p\}$ there exists $g_i\in H$ such that $q(g_i.J_i)\subset I$.
By replacing $J_i$ by $g_i.J_i$ for all $i$ such that $q(J_i)\cap I$ is degenerate, 
we get $q(J_1)\cup\dots \cup q(J_p)=I$.

Let $\phi_i:J_i\ra q(J_i)$ be the partial isometry defined as the restriction of $q$ to $J_i$.
We view $\phi_i$ as a partial isometry whose domain of definition is $\dom\phi_i=J_i$ and 
whose range is $\Im\phi_i=\{1\}\times q(J_i)\subset X(I,S)$.
For $g\in H$, we consider the $g$-conjugate of $\phi_i$ defined by
$$\phi_i^g=g \circ  \phi_i\circ g\m:g.J_i\ra \{g\}\times I.$$
Let $\Phi=\{\phi_i^g|g\in H,i=1,\dots,p\}$.
By construction, $\bigcup_{\phi\in\Phi}\dom\Phi=\Xt$,
and  $\bigcup_{\phi\in\Phi}\Im\Phi=H\times I$.

There is a kind of commutative diagram 
$$\xymatrix@1@R=0.5cm{
\Xt \ar[rr]^-{\Phi}\ar[rd]_q&& H\times I\ar[dl]^p\\
&Y&
}
$$
meaning that for each $x\in \Xt$, and any $\phi\in\Phi$ defined
on $x$, $q(x)=p(\phi(x))$.

We now build a foliated complex $X'$ by gluing foliated bands on $H\times I$.
For each pair of arcs $g_1.J_{i_1}, g_2.J_{i_2}\subset \Xt$ whose intersection $K$ is non-empty,
we add a foliated band joining $\phi_{i_1}^{g_1}(K)$ to
$\phi_{i_2}^{g_2}(K)$ with holonomy  $\phi_{i_2}^{g_2}\circ (\phi_{i_1}^{g_1})\m$.
For each triangle $\tau$ of $X$, and for each pair of arcs 
$g_1.J_{i_1}, g_2.J_{i_2}\subset \partial \tau$ such that the holonomy
along the leaves of $\tau$ defines a partial isometry 
$\psi:K_1\subset g_1.I_1\ra K_2\subset g_2.I_2$,
we add a foliated band 
joining $\phi_{i_1}^{g_1}(K_1)$ to
$\phi_{i_2}^{g_2}(K_2)$, with holonomy
$\phi_{i_2}^{g_2}\circ \psi\circ (\phi_{i_1}^{g_1})\m$.

The important property is that $\Phi$ maps leaves to leaves in the following sense:
if $x,y\in \Xt$ are in the same leaf, and if
$\phi,\psi\in \Phi$ are defined on $x$ and $y$ respectively,
then $\phi(x)$ and $\psi(y)$ are in the same leaf of $X'$.

Since $X$ is locally finite, so is $X'$.
Therefore, the natural free, properly discontinuous action of $H$ on $X'$ is cocompact.
Let $S\subset H$ be the finite set of elements $s\in H$ such that 
there is a band connecting $(1,I)$ to $(s,I)$.
Then $X'$ is naturally contained in $X(I,S)$.
And $\Phi$, viewed as a collection of maps $\Xt\ra X(I,S)$, still maps leaves to leaves.

Let $T(I,S)$ be the leaf space made Hausdorff of $X(I,S)$ and $\Tilde p:X(I,S)\ra T(I,S)$ be the quotient map.
The map $p:H\times I\ra Y$ being constant on any leaf of $X(I,S)$ (or rather on its intersection with $H\times I$), 
and isometric in restriction to any connected component of $H\times I$,
it induces a natural distance decreasing map $f:T(I,S)\ra Y$.
$$\xymatrix@1@R=0.8cm{
\Xt \ar[r]^-{\Phi}\ar[d]_q& H\times I\ar[dl]_p \ar[d]^-{\Tilde p}  \\
Y\ar@/_3mm/@{.>}[r]_-g &T(I,S)\ar[l]_f
}
$$

Since $\Phi$ maps leaves to leaves,
and since $Y$ is the set of leaves of $\Xt$ (it is Hausdorff by Proposition \ref{prop_hausdorff}),
$\Phi$ induces a map $g:Y\ra T(I,S)$.
This map is distance decreasing because 
any path in $X$ written as a concatenation of leaf segments and of transverse arcs in $\Xt$
defines a path having the same transverse measure in $X(I,S)$ as follows:
subdivide its transverse pieces so that they can be mapped to $X(I,S)$ using $\Phi$,
and join the obtained paths by leaf segments of $X(I,S)$.
Since $f$ and $g$ are distance decreasing,
both are isometries. 
In particular $Y$ is dual to $X(I,S)$.

We prove that the leaf space of $X(I,S)$ is Hausdorff.
Assume that $x,y\in H\times I$ have the image under $p$.
Consider $\Tilde x,\Tilde y\in X$ and $\phi,\psi\in\Phi$ such that
$\phi(\Tilde x)=x$ and $\psi(\Tilde y)=y$.
Then $q(\Tilde x)=q(\Tilde y)$, so $\Tilde x$ and $\Tilde y$ 
are in the same leaf because the leaf space of $X$ is Hausdorff.
Therefore, $x$ and $y$ are in the same leaf of $X(I,S)$.
\end{proof}

\begin{rem*}
  If the lemma holds for some $S\subset H$, it also holds for any $S'$ containing $S$.
\end{rem*}

\subsection{Pseudo-groups}

Up to now we considered partial isometries between closed intervals.
We now need to consider partial isometries between \emph{open} intervals.
We use notations like $\rond D,\rond \phi$ to emphasize this point.
An open \emph{multi-interval} $\rond D$ is a finite union of copies of bounded open intervals of $\bbR$.
The \emph{pseudo-group of isometries} generated by some partial isometries $\rond \phi_1,\dots,\rond \phi_n$
between open intervals of $D$,
is the set of partial isometries $\rond\phi:\rond I\ra\rond J$ such that for any $x\in \rond I$, there exists
a composition $\rond\phi_{i_1}^{\pm 1}\circ\dots\circ \rond\phi_{i_k}^{\pm 1}$
which is defined on a neighbourhood of $x$, and which coincides with $\rond\phi$ on this neighbourhood.
Two points $x,y\in\rond{D}$ are in the same orbit if there exists $\rond\phi\in\Lambda$ such that $y=\rond\phi(x)$.
A pseudo-group is \emph{minimal} if its orbits are dense in $\rond D$.

Consider an arc $I\subset Y$, where $H\actson Y$ is an indecomposable nice action.
Let $\rond I=I\setminus \partial I$ and
 $\rond\phi_s$ be the restriction of $s$ to  the interior  of $s\m.\rond I\cap \rond I$.
Let $\Lambda_0$ be the pseudo-group of isometries on $\rond{I}$
generated $\{\rond\phi_s|s\in S\}$.
Let $\rond D\subset \rond{I}$ be the set of points whose orbit under $\Lambda_0$ is infinite.

Let $\Lambda(I,S)$ be the restriction of $\Lambda_0$ to $\rond D$.
We claim that $\rond I\setminus \rond D$ is finite (in particular, $\rond D$ is an open multi-interval)
and that $\Lambda(I,S)$ is minimal.
Let $C$ be the cut locus of $\Sigma=X(I,S)/H$.
Since $Y$ is indecomposable, $\Sigma\setminus C$ consists of only one minimal component,
so every leaf $l$ of $\Sigma\setminus C$ is dense in $\Sigma\setminus C$.
Since for any such leaf $l$, $l\cap \rond I$ is contained in a $\Lambda_0$-orbit,
the claim follows.

A pseudo-group of isometries $\Lambda$ is \emph{homogeneous}
if for each $\rond\phi\in\Lambda$, any partial isometry $\rond\psi$ extending $\rond\phi$ lies in
$\Lambda$. 
When $\Lambda(I,S)$ is homogeneous, $Y$ is a line (see for instance \cite[section 5]{Gui_approximation}).

\begin{thm}[\cite{Lev_pseudo,Gui_these}]\label{thm_minimax}
  Let $\Lambda$ be a minimal finitely generated pseudo-group of isometries of an open multi-interval $\rond D$.
Then the set of non-homogeneous pseudo-groups of isometries containing $\Lambda$ is finite.
\end{thm}

\begin{rem*}
  In \cite{Lev_pseudo}, the result is proved for orientable pseudogroups of isometries of the circle.
Theorem \ref{thm_minimax} is an easy generalisation, proved in \cite{Gui_these}, but is not published.
The proof follows step by step the proof in \cite{Lev_pseudo}, 
using Gusm\~ao's extension of Levitt's results (\cite{Gus_feuilletages}).
If $\Lambda$ is non-orientable, one can deduce Theorem \ref{thm_minimax} from the orientable case
by a straightforward 2-fold covering argument. The only remaining unpublished case is when $\Lambda$ is orientable, and we allow
larger pseudo-groups to be non-orientable, but this case is not needed in our argument.
\end{rem*}

\subsection{Stabilisation of indecomposable components}

We are now ready to prove the stabilisation of indecomposable components.

\begin{proof}[Proof of Lemma \ref{lem_inj}]
Let $G_k\actson T_k$ be a sequence of geometric actions converging strongly to $G\actson T$
as in lemma \ref{lem_cv2}. 
Recall that $(\phi_k,f_k):G_k\actson T_k\ra G_{k+1}\actson T_{k+1}$ and $(\Phi_k,F_k):G_k\actson T_k\ra G\actson T$
denote the maps of the corresponding direct system.

Let $H_{k_0}\actson Y_{k_0}$ be an indecomposable component of $T_{k_0}$.
Consider $k\geq k_0$.
Since $f_{{k_0}k}(Y_{k_0})$ is indecomposable, it is contained in an indecomposable component $Y_k$
of $T_k$ (see Lemma \ref{lem_indec_component} and \ref{lem_indecompo}).
Let $I_{k_0}\subset Y_{k_0}$ be an arc which embeds into $T$ under $F_{k_0}$,
and let $I_k=f_{{k_0}k}(I_{k_0})$.

By Lemma \ref{lem_dual}, $Y_k$ is dual to $X(I_k,S_k)$ for some finite set $S_k\subset H_k$.
By enlarging each $S_k$, we can assume that for all $k$, $\phi_k(S_k)\subset S_{k+1}$.
Let $\Lambda(I_k,S_k)$ be the corresponding minimal pseudo-group of isometries on $\rond D_k\subset \rond I_k$.
Under the natural identification between $I_k$ and $I_{k+1}$, we get
$\rond D_{k}\subset \rond D_{k+1}$, and  $\Lambda(I_{k},S_{k})\subset \Lambda(I_{k+1},S_{k+1})$.
Since $I_k\setminus \rond D_k$ is finite, for $k$ large enough, $\rond D_{k+1}=\rond D_k$.

By Theorem \ref{thm_minimax} above, for $k$ large enough, either $\Lambda(I_k,S_k)$ is homogeneous,
or $\Lambda(I_k,S_k)=\Lambda(I_{k+1},S_{k+1})$.
It is an exercise to show that if $\Lambda$ is a minimal homogeneous pseudo-group of isometries,
then any pseudo-group of isometries containing $\Lambda$ is homogeneous.
Thus, in the first case, $\Lambda(I_{k+1},S_{k+1})$ is homogeneous.
It follows that $Y_k$ and $Y_{k+1}$ are lines. Since $H_k$ acts with dense orbits on $Y_k$,
the morphism of $\bbR$-trees $f_{k}{}_{|Y_k}:Y_k\ra Y_{k+1}$ 
is necessarily one-to-one: there exists an arc of $Y_k$ which is embedded under $f_{k}$ and using the action of $H_k$,
we see that $f_k$ is locally isometric. It follows that $f_k(Y_k)=Y_{k+1}$.

In the second case, the following result concludes the proof.
\end{proof}

\begin{lem}[The pseudo-group determines the action]
Assume that $(\phi,f):H\actson Y\ra H'\actson Y'$ 
is a morphism between nice indecomposable actions.
Consider $I\subset Y$ such that $f$ is isometric in restriction to $I$ and let $I'=f(I)$.
Consider finite subsets $S\subset H$ and $S'\subset H'$ such that  $\phi(S)\subset S'$.

Assume that $Y$ is dual to $X(I,S)$, that $Y'$ is dual to $X(I',S')$, 
and that $\Lambda(I,S)=\Lambda(I',S')$.

Then $f:Y\ra Y'$ is a surjective isometry and $\phi:H\ra H'$ is an isomorphism.
\end{lem}

\begin{figure}[htb]
  \centering
  \includegraphics{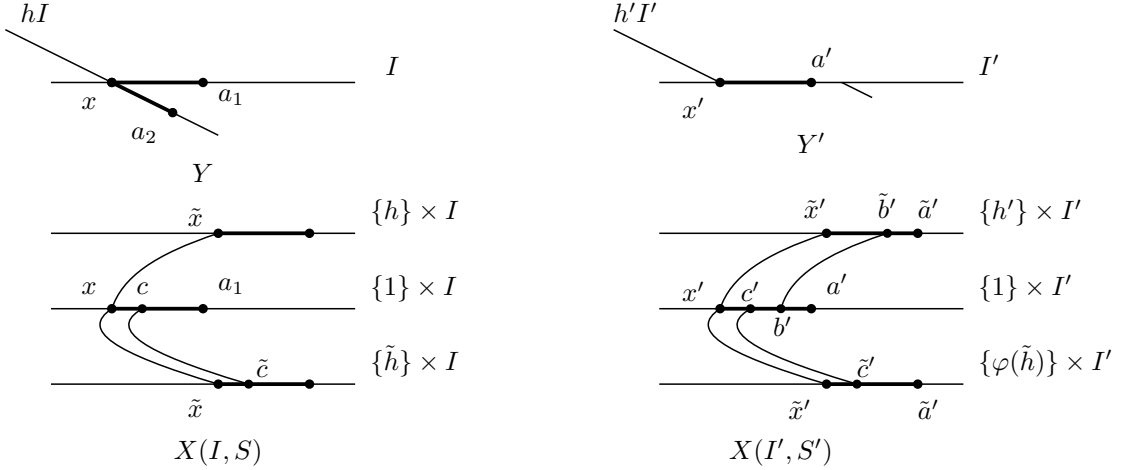}
  \caption{The pseudo-group determines the action}
  \label{fig_121}
\end{figure}
 
\begin{proof}
 Assume that $f$ is not one-to-one.
There exists two arcs $[x,a_1],[x,a_2]\subset Y$ such that $[x,a_1]\cap [x,a_2]=\{x\}$ and $f([x,a_1])=f([x,a_2])$.
Since $I$ spans $Y$, one may shorten these arcs and assume without loss of generality that
$[x,a_1]\subset I$ and $[x,a_2]\subset h.I$ for some $h\in H$.

Consider $a'=f(a_1)=f(a_2)$, $x'=f(x)$ and $h'=\phi(h)$ (see Figure \ref{fig_121}).
The two points $\Tilde a'=h'{}\m a'$ and $\Tilde x'=h'{}\m x'$ lie in $I'$.
The subsets $\{1\}\times [x',a']$ and $\{h'\}\times \Tilde [\tilde x',\Tilde a']$ of $X(I',S')$ map to the same arc in $Y'$.
By Corollary \ref{cor_bands},
 there exist arcs $[x',b']\subset [x',a']$, $[\Tilde x',\Tilde b']\subset [\Tilde x',\Tilde a']$,
and a holonomy band in $X(I',S')$ joining $\{1\}\times [x',b']$
to $\{h'\}\times [\Tilde x',\Tilde b']$.
This holonomy band defines a partial isometry $\rond\psi{}'\in\Lambda(I',S')$
such that $\rond\psi{}'((x',b'))=(\Tilde x',\Tilde b')$.
By hypothesis, $\rond\psi{}'$ corresponds to an element of $\rond\psi\in\Lambda(I,S)$,
sending $(x,b)$ to $(\Tilde x,\Tilde b)$ where $b,\Tilde x,\Tilde b$ are the points of $I$
corresponding to $b',\Tilde x',\Tilde b'$ under the natural identification.

We now apply segment closed property to the partial isometry induced by $\rond\psi$ in $X(I,S)/H$.
We get the existence of arcs $[x,c]\subset [x,b]$
and $[\Tilde x,\Tilde c]\subset [\Tilde x,\Tilde b]$
and of a holonomy band in $X(I,S)/H$ whose holonomy
coincides with $\rond\psi_{|(x,c)}$.
Lifting this holonomy band to $X(I,S)$, we get a holonomy band joining
$\{1\}\times [x,c]$ to $\{\Tilde h\}\times [\Tilde x,\Tilde c]$ for some $\Tilde h\in H$.
Since $S'\supset \phi(S)$,
this holonomy band defines a holonomy band in $X(I',S')$
joining $\{1\}\times [x',c']$ to $\{\phi(\Tilde h)\}\times [\Tilde x',\Tilde c']$
where $\Tilde c'=f(\Tilde c)$.
In particular, $\phi(\Tilde hh\m)$ fixes the arc $[\Tilde x',\Tilde c']$, so $\phi(\Tilde hh\m)=1$.
But since $\{\Tilde h\}\times\{ \Tilde x\}$ is in the same leaf as $\{h\}\times \{\Tilde x\}$,
$\Tilde hh\m$ fixes the point $x$ in $Y$.
Since $\phi$ is one-to-one in restriction to point stabilizers of $Y$,
$h=\Tilde h$. Since $\{h\}\times [\Tilde x,\Tilde c]$ maps into $[x,a_2]$ in $Y$,
we get that $[x,a_2]\cap [x,a_1]$ contains $[x,c]$, a contradiction.

This proves that $f$ is one-to-one.
It follows that $\phi$ is one-to one since an element of its kernel
must fix $Y$ pointwise.

Let's prove that $\phi(H)\supset H'$, the other inclusion being trivial.
The fact that $f(Y)=Y'$ will follow because $H'\actson Y'$  is minimal. 

Since $H'\actson Y$ is indecomposable, $H'$ is generated by the set of elements $g\in H'$ such that $g.I'\cap I'$ is
non-degenerate.
Let $h'\in H'$ be such that $[x',a']=h'.I'\cap I'$ is non-degenerate.
Repeating the argument above, there exists 
an arc $[x',b']\subset [x',a']$ and
$\rond\psi{}'\in\Lambda(I',S')$, 
$\rond\psi{}':(x',b')\ra (\Tilde x',\Tilde b')$ where $\Tilde x'=h\m x'$, $\Tilde b'=h\m b'$.
Let $\rond\psi\in\Lambda(I,S)$ be the corresponding partial isometry.
By segment closed property, there exist
arcs $[x,c]\subset [x,b]$ and $[\Tilde x,\Tilde c]\subset[\Tilde x,\Tilde b]$
and a holonomy band in $X(I,S)/H$  joining them in $X(I,S)/H$, 
whose holonomy coincides with $\rond\psi$.
This band lifts in $X(I,S)$, to a band  joining 
$\{1\}\times [x,c]$ to $\{\Tilde h\}\times [\Tilde x,\Tilde c]$ for some $\Tilde h\in H$.
In $X(I',S')$, this bands joins the corresponding arcs 
$\{1\}\times [x,c]$ to $\{\Tilde h\}\times [\Tilde x,\Tilde c]$.
It follows that in $Y$, the actions of $\phi(\Tilde h)$ and $h$ coincide on $[x',c']$.
Since arc stabilizers are trivial, $h=\phi(\Tilde h)$.
We conclude that $\phi(H)=H'$.
\end{proof}

\small

\bibliographystyle{alpha}
\bibliography{published,unpublished}

\def\cprime{$'$}
\begin{thebibliography}{Gui98b}

\bibitem[Ali04]{Alibegovic_MR}
Emina Alibegovic.
\newblock Makanin-razborov diagrams for limit groups.
\newblock math.GR/0410198, 2004.

\bibitem[BF]{BF_outer}
Mladen Bestvina and Mark Feighn.
\newblock Outer limits.
\newblock Preprint.

\bibitem[BF91]{BF_complexity}
Mladen Bestvina and Mark Feighn.
\newblock Bounding the complexity of simplicial group actions on trees.
\newblock {\em Invent. Math.}, 103(3):449--469, 1991.

\bibitem[BF92]{BF_combination}
M.~Bestvina and M.~Feighn.
\newblock A combination theorem for negatively curved groups.
\newblock {\em J. Differential Geom.}, 35(1):85--101, 1992.

\bibitem[BF95]{BF_stable}
Mladen Bestvina and Mark Feighn.
\newblock Stable actions of groups on real trees.
\newblock {\em Invent. Math.}, 121(2):287--321, 1995.

\bibitem[BS05]{BeleSzcz_endomorphisms}
Igor Belegradek and Andrzej Szczepanski.
\newblock Endomorphisms of relatively hyperbolic groups.
\newblock math.GR/0501321, 2005.

\bibitem[Chi01]{Chi_book}
Ian Chiswell.
\newblock {\em Introduction to $\Lambda$-trees}.
\newblock World Scientific Publishing Co. Inc., River Edge, NJ, 2001.

\bibitem[CL95]{CoLu_very}
Marshall~M. Cohen and Martin Lustig.
\newblock Very small group actions on {$\mathbb{R}$}-trees and {D}ehn twist
  automorphisms.
\newblock {\em Topology}, 34(3):575--617, 1995.

\bibitem[CM87]{CuMo}
Marc Culler and John~W. Morgan.
\newblock Group actions on {$\mathbb{R}$}-trees.
\newblock {\em Proc. London Math. Soc. (3)}, 55(3):571--604, 1987.

\bibitem[Coo05]{Magnus_software}
Group~Theory Cooperative.
\newblock Magnus, computational package for exploring infinite groups, version
  4.1.3 beta, June 2005.
\newblock (G.Baumslag director).

\bibitem[Del99]{Delzant_accessibilite}
Thomas Delzant.
\newblock Sur l'accessibilit\'e acylindrique des groupes de pr\'esentation
  finie.
\newblock {\em Ann. Inst. Fourier (Grenoble)}, 49(4):1215--1224, 1999.

\bibitem[DS06]{DrutuSapir_tree-graded}
Cornelia Dru{\c t}u and Mark Sapir.
\newblock Groups acting on tree-graded spaces and splittings of relatively
  hyperbolic group.
\newblock math.GR/0601305, 2006.

\bibitem[Dun98]{Dun_folding}
M.~J. Dunwoody.
\newblock Folding sequences.
\newblock In {\em The Epstein birthday schrift}, pages 139--158 (electronic).
  Geom. Topol., Coventry, 1998.

\bibitem[GLP94]{GLP1}
D.~Gaboriau, G.~Levitt, and F.~Paulin.
\newblock Pseudogroups of isometries of {$\mathbb{R}$} and {R}ips' theorem on
  free actions on {$\mathbb{R}$}-trees.
\newblock {\em Israel J. Math.}, 87(1-3):403--428, 1994.

\bibitem[GLP95]{GLP2}
Damien Gaboriau, Gilbert Levitt, and Fr{\'e}d{\'e}ric Paulin.
\newblock Pseudogroups of isometries of {$\mathbb{R}$}: reconstruction of free
  actions on {$\mathbb{R}$}-trees.
\newblock {\em Ergodic Theory Dynam. Systems}, 15(4):633--652, 1995.

\bibitem[Gro05]{Groves_limit_hypII}
Daniel Groves.
\newblock Limit groups for relatively hyperbolic groups. {II}.
  {M}akanin-{R}azborov diagrams.
\newblock {\em Geom. Topol.}, 9:2319--2358 (electronic), 2005.

\bibitem[Gui98a]{Gui_these}
Vincent Guirardel.
\newblock {\em Actions de groupes sur des arbres r{\'e}els et dynamique dans la
  fronti{\`e}re de l'outre-espace}.
\newblock PhD thesis, Universit{\'e} {T}oulouse III, January 1998.

\bibitem[Gui98b]{Gui_approximation}
Vincent Guirardel.
\newblock Approximations of stable actions on {$\mathbb{R}$}-trees.
\newblock {\em Comment. Math. Helv.}, 73(1):89--121, 1998.

\bibitem[Gui04]{Gui_limit}
Vincent Guirardel.
\newblock Limit groups and groups acting freely on {$\Bbb R\sp n$}-trees.
\newblock {\em Geom. Topol.}, 8:1427--1470 (electronic), 2004.

\bibitem[Gui05]{Gui_coeur}
Vincent Guirardel.
\newblock C\oe ur et nombre d'intersection pour les actions de groupes sur les
  arbres.
\newblock {\em Ann. Sci. \'Ecole Norm. Sup. (4)}, 38(6):847--888, 2005.

\bibitem[Gus00]{Gus_feuilletages}
Paulo Gusm{\~a}o.
\newblock Feuilletages mesur\'es et pseudogroupes d'isom\'etries du cercle.
\newblock {\em J. Math. Sci. Univ. Tokyo}, 7(3):487--508, 2000.

\bibitem[Ima79]{Imanishi}
Hideki Imanishi.
\newblock On codimension one foliations defined by closed one-forms with
  singularities.
\newblock {\em J. Math. Kyoto Univ.}, 19(2):285--291, 1979.

\bibitem[KW05]{KaWe_acylindrical}
Ilya Kapovich and Richard Weidmann.
\newblock Acylindrical accessibility for groups acting on {$\Bbb R$}-trees.
\newblock {\em Math. Z.}, 249(4):773--782, 2005.

\bibitem[Lev93]{Lev_pseudo}
Gilbert Levitt.
\newblock La dynamique des pseudogroupes de rotations.
\newblock {\em Invent. Math.}, 113(3):633--670, 1993.

\bibitem[Lev94]{Lev_graphs}
Gilbert Levitt.
\newblock Graphs of actions on {$\mathbb {R}$}-trees.
\newblock {\em Comment. Math. Helv.}, 69(1):28--38, 1994.

\bibitem[LP97]{LP}
Gilbert Levitt and Fr{\'e}d{\'e}ric Paulin.
\newblock Geometric group actions on trees.
\newblock {\em Amer. J. Math.}, 119(1):83--102, 1997.

\bibitem[Mor88]{Mo_ergodic}
John~W. Morgan.
\newblock Ergodic theory and free actions of groups on {$\mathbb {R}$}-trees.
\newblock {\em Invent. Math.}, 94(3):605--622, 1988.

\bibitem[MS84]{MS_valuationsI}
John~W. Morgan and Peter~B. Shalen.
\newblock Valuations, trees, and degenerations of hyperbolic structures. {I}.
\newblock {\em Ann. of Math. (2)}, 120(3):401--476, 1984.

\bibitem[Pau88]{Pau_topologie}
Fr{\'e}d{\'e}ric Paulin.
\newblock Topologie de {G}romov \'equivariante, structures hyperboliques et
  arbres r\'eels.
\newblock {\em Invent. Math.}, 94(1):53--80, 1988.

\bibitem[RS94]{RiSe_structure}
E.~Rips and Z.~Sela.
\newblock Structure and rigidity in hyperbolic groups. {I}.
\newblock {\em Geom. Funct. Anal.}, 4(3):337--371, 1994.

\bibitem[Sco73]{Scott_coherent}
G.~P. Scott.
\newblock Finitely generated $3$-manifold groups are finitely presented.
\newblock {\em J. London Math. Soc. (2)}, 6:437--440, 1973.

\bibitem[Sel97]{Sela_acylindrical}
Z.~Sela.
\newblock Acylindrical accessibility for groups.
\newblock {\em Invent. Math.}, 129(3):527--565, 1997.

\bibitem[Sel99]{Sela_hopf}
Z.~Sela.
\newblock Endomorphisms of hyperbolic groups. {I}. {T}he {H}opf property.
\newblock {\em Topology}, 38(2):301--321, 1999.

\bibitem[Sel01]{Sela_diophantine1}
Zlil Sela.
\newblock Diophantine geometry over groups. {I}. {M}akanin-{R}azborov diagrams.
\newblock {\em Publ. Math. Inst. Hautes \'Etudes Sci.}, 93:31--105, 2001.

\bibitem[Sel02]{Sela_diophantine7}
Zlil Sela.
\newblock Diophantine geometry over groups vii: The elementary theory of a
  hyperbolic group.
\newblock {\url{http://www.ma.huji.ac.il/~zlil}}, 2002.

\bibitem[Sel06]{Sela_diophantine6}
Z.~Sela.
\newblock Diophantine geometry over groups. {VI}. {T}he elementary theory of a
  free group.
\newblock {\em Geom. Funct. Anal.}, 16(3):707--730, 2006.

\bibitem[Ser77]{Serre_arbres}
Jean-Pierre Serre.
\newblock {\em Arbres, amalgames, ${\rm {S}{L}}\sb{2}$}.
\newblock Soci\'et\'e Math\'ematique de France, Paris, 1977.
\newblock R\'edig\'e avec la collaboration de Hyman Bass, Ast\'erisque, No. 46.

\bibitem[Sha91]{Sh_dendrology}
Peter~B. Shalen.
\newblock Dendrology and its applications.
\newblock In {\em Group theory from a geometrical viewpoint (Trieste, 1990)},
  pages 543--616. World Sci. Publishing, River Edge, NJ, 1991.

\bibitem[Sko89]{Skora_combination}
Richard Skora.
\newblock Combination theorems for actions on trees.
\newblock preprint, 1989.

\bibitem[Sta83]{Stallings_topology}
John~R. Stallings.
\newblock Topology of finite graphs.
\newblock {\em Invent. Math.}, 71(3):551--565, 1983.

\bibitem[Swa04]{Swarup_Delzant}
Gadde Swarup.
\newblock {D}elzant's variation on {S}cott complexity.
\newblock arXiv:math.GR/0401308, 2004.

\end{thebibliography}

\begin{flushleft}
Vincent Guirardel\\
Institut de Mathématiques de Toulouse, UMR 5219\\
Universit\'e Paul Sabatier\\
31062 Toulouse cedex 9.\\
France.\\
\emph{e-mail:}\texttt{guirardel@math.ups-tlse.fr}\\
\end{flushleft}

\end{document}